\documentclass[11pt]{article}

\usepackage[a4paper,margin=2.5cm]{geometry}
\usepackage[T1]{fontenc}
\usepackage{amsmath}
\usepackage{amsthm}
\usepackage{newtxtext}
\usepackage{newtxmath}
\usepackage{graphicx}
\usepackage[round,authoryear]{natbib}
\usepackage{tikz}
\usetikzlibrary{plotmarks}
\usepackage[hidelinks]{hyperref}
\hypersetup{pdftitle={An iterated-logarithm price of anytime validity for log-optimal e-processes},pdfauthor={Adrien Corenflos}}

\newtheorem{theorem}{Theorem}
\newtheorem{proposition}{Proposition}
\newtheorem{lemma}{Lemma}

\theoremstyle{definition}

\newtheorem{example}{Example}
\theoremstyle{remark}
\newtheorem{remark}{Remark}

\makeatletter
\renewenvironment{proof}[1][]{\par\pushQED{\qed}\normalfont\topsep6\p@\@plus6\p@\relax\trivlist
  \item[\hskip\labelsep\itshape\proofname\if\relax\detokenize{#1}\relax\else\ #1\fi\@addpunct{.}]\ignorespaces}
  {\popQED\endtrivlist\@endpefalse}
\makeatother
\newcommand{\tbl}[2]{\caption{#1}\centering#2}

\newcommand{\E}{\mathbb{E}}
\newcommand{\KL}{\mathrm{KL}}
\newcommand{\N}{\mathrm{N}}
\newcommand{\cP}{\mathcal{P}}
\newcommand{\cF}{\mathcal{F}}
\newcommand{\cG}{\mathcal{G}}
\newcommand{\cB}{\mathcal{B}}
\newcommand{\cE}{\mathcal{E}}
\newcommand{\cX}{\mathcal{X}}
\newcommand{\R}{\mathbb{R}}
\newcommand{\Es}{E^{\star}}
\newcommand{\Ps}{P^{\star}}
\newcommand{\dd}{\mathrm{d}}
\newcommand{\ind}{\mathbf{1}}
\newcommand{\pr}{\mathbb{P}}
\makeatletter
\newenvironment{parts}{\if@inlabel\leavevmode\par\fi\protected@edef\parts@outer{\@currentlabel}\renewcommand\theenumi{(\roman{enumi})}\renewcommand\p@enumi{\parts@outer}\@ifundefined{enummax}{\enumerate}{\enumerate[(iii)]}}{\endenumerate}
\newcommand\partlabel[1]{\label{#1}{\protected@edef\@currentlabel{\theenumi}\label{#1*}}}
\makeatother
\DeclareRobustCommand{\legEtaThree}{\tikz[baseline=-0.6ex]{\draw[color=black, line width=0.8pt, solid] (0,0) -- (16pt,0); \draw[color=black, line width=0.6pt, fill=black] plot[mark=*, mark size=1.6pt] coordinates {(8pt,0)};}}
\DeclareRobustCommand{\legLowThree}{\tikz[baseline=-0.6ex]{\draw[color=black, line width=0.8pt, dashed] (0,0) -- (16pt,0); \draw[color=black, line width=0.6pt, fill=white] plot[mark=*, mark size=1.6pt] coordinates {(8pt,0)};}}
\DeclareRobustCommand{\legLimThree}{\tikz[baseline=-0.6ex]{\draw[color=black, line width=0.8pt, dotted] (0,0) -- (16pt,0);}}
\DeclareRobustCommand{\legEtaFive}{\tikz[baseline=-0.6ex]{\draw[color={rgb,255:red,85;green,85;blue,85}, line width=0.8pt, solid] (0,0) -- (16pt,0); \draw[color={rgb,255:red,85;green,85;blue,85}, line width=0.6pt, fill={rgb,255:red,85;green,85;blue,85}] plot[mark=square*, mark size=1.6pt] coordinates {(8pt,0)};}}
\DeclareRobustCommand{\legLowFive}{\tikz[baseline=-0.6ex]{\draw[color={rgb,255:red,85;green,85;blue,85}, line width=0.8pt, dashed] (0,0) -- (16pt,0); \draw[color={rgb,255:red,85;green,85;blue,85}, line width=0.6pt, fill=white] plot[mark=square*, mark size=1.6pt] coordinates {(8pt,0)};}}
\DeclareRobustCommand{\legLimFive}{\tikz[baseline=-0.6ex]{\draw[color={rgb,255:red,85;green,85;blue,85}, line width=0.8pt, dotted] (0,0) -- (16pt,0);}}
\DeclareRobustCommand{\legEtaSeven}{\tikz[baseline=-0.6ex]{\draw[color={rgb,255:red,153;green,153;blue,153}, line width=0.8pt, solid] (0,0) -- (16pt,0); \draw[color={rgb,255:red,153;green,153;blue,153}, line width=0.6pt, fill={rgb,255:red,153;green,153;blue,153}] plot[mark=triangle*, mark size=1.6pt] coordinates {(8pt,0)};}}
\DeclareRobustCommand{\legLowSeven}{\tikz[baseline=-0.6ex]{\draw[color={rgb,255:red,153;green,153;blue,153}, line width=0.8pt, dashed] (0,0) -- (16pt,0); \draw[color={rgb,255:red,153;green,153;blue,153}, line width=0.6pt, fill=white] plot[mark=triangle*, mark size=1.6pt] coordinates {(8pt,0)};}}
\DeclareRobustCommand{\legLimSeven}{\tikz[baseline=-0.6ex]{\draw[color={rgb,255:red,153;green,153;blue,153}, line width=0.8pt, dotted] (0,0) -- (16pt,0);}}
\title{\textbf{An iterated-logarithm price of anytime validity for log-optimal e-processes}}
\author{Adrien Corenflos\\{\small Department of Statistics, University of Warwick, Coventry CV4 7AL, U.K.}\\{\small\texttt{adrien.corenflos@warwick.ac.uk}}}
\date{}

\begin{document}

\maketitle

\begin{abstract}
Log-optimal e-variables at fixed horizons, also called numéraires, have been conjectured to ``almost'' form an e-process for exponential-family null hypotheses.
  A related open question is whether an e-process can attain the optimal e-power at every horizon up to a bounded loss.
  Given a final horizon, we study the minimax regret of e-processes relative to the numéraires at all times up to it.
  For the Gaussian location null and the composite alternative of Gaussian distributions with a smaller variance and an unknown location, we show that the minimax regret in the worst case over the location grows like the iterated logarithm of the final horizon.
  A single e-process attains this rate at every horizon.
  This growth answers the question negatively.
  It also refutes the conjecture if ``almost'' means a regret bounded uniformly over horizons.
  For this Gaussian location problem, we thus exhibit the frontier precisely: the minimax regret is bounded if and only if the alternatives have at least the variance of the null.

\end{abstract}

\section{Introduction}\label{sec:intro}

\subsection{E-processes and the minimax regret}\label{sec:intro-regret}

Let $X_1,X_2,\ldots$ be observations with values in a measurable space $\cX$, let $X^t=(X_1,\ldots,X_t)$, and let $(\cF_t)_{t\ge0}$ be the filtration generated by $X=(X_t)_{t\ge1}$, with $\cF_0$ trivial.
Let $\cP$ be a set of null distributions and $Q$ an alternative distribution of $X$.
We define e-variables and e-processes by the inequalities
\begin{equation}\label{eq:eproc}
  \smash[t]{
    \sup_{P\in\cP}\E_P[E]\le1,\quad\sup_{P\in\cP}\E_P[E_\tau]\le1,
  }
\end{equation}
where $\E_P$ is the expectation under $P$, $\tau$ ranges over the stopping times of $(\cF_t)$ with values in $\{0,1,\ldots\}\cup\{\infty\}$, and $E_\infty=\limsup_{t\to\infty}E_t$.
An e-variable at horizon $t$ is an $\cF_t$-measurable random variable $E$ with values in $[0,\infty]$ that satisfies the first inequality in~\eqref{eq:eproc}~\citep{vovk2021evalues}, while an e-process is a process $(E_t)_{t\ge0}$ adapted to $(\cF_t)$ with values in $[0,\infty]$ that satisfies the second \citetext{\citealp{grunwald2024safe}; \citealp[Section~2.4]{ramdas2023game}; \citealp[Definition~7.3]{ramdas2025hypothesis}}.
By Markov's inequality, we can reject the null at level $\alpha$ as soon as an e-process reaches $1/\alpha$~\citep[Sections~2.5 and~2.6]{ramdas2023game}.

Since the constant $E_t \equiv 1$ is an e-process, we also measure the power of e-variables under the alternative $Q$.
Let $\cE_t$ be the set of e-variables at horizon $t$, and let $P_t$ be the distribution of $X^t$ under a distribution $P$ of $X$.
The e-power of $E\in\cE_t$ is $\E_Q[\log E]$.
To describe the e-variables with the largest e-power, \citet[Definition~2.1]{larsson2025numeraire} call $\Es_t\in\cE_t$ a numéraire at horizon $t$ if it satisfies the first of the properties
\begin{equation}\label{eq:numeraire}
  \sup_{E\in\cE_t}\E_Q[E/\Es_t]\le1,\quad\E_Q[\log\Es_t]=\sup_{E\in\cE_t}\E_Q[\log E],\quad\Es_t=\dd Q_t/\dd\Ps_t.
\end{equation}
They prove that a numéraire exists at each horizon for every null and every alternative, and that it satisfies the second property in~\eqref{eq:numeraire}~\citep[Theorem~2.6 and Proposition~2.4]{larsson2025numeraire}.
When $Q_t$ is absolutely continuous with respect to the null, the numéraire also satisfies the third property, for a sub-probability measure $\Ps_t$ on $\cX^t$~\citep[Theorem~2.6 and Lemma~3.3]{larsson2025numeraire}.
When the null is composite, $\Ps_t$ may be a mixture of null distributions with a mixing distribution that changes with $t$, as in Proposition~\ref{prop:numeraire}.
The sequence $(\Es_t)_{t\ge0}$ is then in general not an e-process.
We therefore compare e-processes with the numéraires through the regret of an e-variable $E\in\cE_t$ at $t$ and the minimax regret up to a horizon $T$,
\begin{equation}\label{eq:price}
  R_t(E)=\E_Q[\log(\Es_t/E)]\ge0,\quad\Gamma_T=\inf_{(E_t)}\max_{0\le t\le T}R_t(E_t),
\end{equation}
where the infimum runs over all e-processes $(E_t)$, and we call $R_t(E_t)$ the regret of the e-process $(E_t)$ at $t$.

\subsection{Known results and open questions}\label{sec:intro-known}

We first recall why $\Gamma_T=0$ for a simple null.
A test (super)martingale with respect to a filtration is a nonnegative process with initial value one that is a (super)martingale with respect to this filtration under every null distribution.
Test supermartingales with respect to $(\cF_t)$ are e-processes.
For a simple null $\cP=\{P\}$, the numéraires $\Es_t=\dd Q_t/\dd P_t$ form a test martingale with respect to $(\cF_t)$.
Consequently, $\Gamma_T=0$ for every $T$~\citep[Theorem~8]{koolen2022log}.
In contrast, for a composite null, we observe in Proposition~\ref{prop:basic:i} that $\inf_{(E_t)}R_t(E_t)=0$ holds at each single horizon $t$ separately.
Similarly, \citet{koning2026anytime} state that ``anytime validity is free'', because a valid test at a planned horizon $T$ induces an anytime-valid sequential test that coincides with it at $T$~\citep[Theorems~1 and~2 and Proposition~3]{koning2026anytime}.
Like Proposition~\ref{prop:basic:i}, their statement holds at a single planned horizon.
Moreover, over all horizons up to $T$, the average $(T+1)^{-1}\sum_{i=0}^{\min(t,T)}\Es_i$ of the numéraires is an e-process~\citep[Section~7.9, Proposition~7.23]{ramdas2025hypothesis}.
We show in Proposition~\ref{prop:basic:ii} that its regret is at most $\log(T+1)$ at every horizon $t\le T$.
However, it is unknown whether some e-process has a regret bounded uniformly over all horizons.
\citet[Section~5]{grunwald2024safe} ask under what conditions a test supermartingale with respect to $(\cF_t)$ has such a bounded regret.
\citet[Sections~1 and~7 and Appendix~E]{hao2026evalues} note that their conditional e-variables usually do not define an e-process, but they nevertheless conjecture informally that these e-variables ``almost'' form one.

Asymptotic growth-rate optimality depends only on the limit of $t^{-1}\log E_t$ \citetext{\citealp[Definition~2.1]{waudbysmith2025instance}; \citealp[Theorem~3.5]{ram2026growth}}.
Since two e-processes with regrets differing by $o(t)$ have the same limit in expectation, this criterion cannot distinguish a bounded regret from one that grows like $\log\log T$.
To our knowledge, no lower bound on $\Gamma_T$ growing with $T$ exists.

\subsection{Contributions}\label{sec:intro-contributions}

We give two general bounds on $\Gamma_T$, valid for every null and every alternative.
From above, $\Gamma_T\le\log\eta_T$, where $\eta_T$ is the largest expected value, under a null distribution, of a numéraire stopped at a stopping time bounded by $T$ (Proposition~\ref{prop:basic:iii}).
From below, $\Gamma_T$ is controlled by how closely mixtures of null distributions stopped at random times can approximate the alternative stopped at a random horizon \citetext{Lemma~\ref{lem:duality}; \citealp[Theorem~5]{koolen2022log}}.
Such mixtures give a positive lower bound only if their stopping times depend on the data (Remark~\ref{rem:deterministic}).

Our main result is concerned with the Gaussian location null $\{\N(\theta,1)^{\otimes\infty}:\theta\in\R\}$ and the composite alternative $\{\N(\mu,s^2)^{\otimes\infty}:\mu\in\R\}$ (Section~\ref{sec:gauss}).
Here $\N(m,v)$ is the normal distribution with mean $m$ and variance $v$, and $P^{\otimes\infty}$ is the distribution of a sequence of independent variables with distribution $P$.
We measure the regret~\eqref{eq:price} against each alternative $\N(\mu,s^2)^{\otimes\infty}$ and take the worst case over $\mu$.
For a known location, the minimax regret over all horizons is zero if $s\le1$ and at most $\log s$ otherwise (Proposition~\ref{prop:simple}).
An unbounded minimax regret can therefore come only from the unknown location.
If $s<1$, the minimax regret at a single horizon is $\kappa=(s^2-1-\log s^2)/2$.
The likelihood ratio of the maximal invariant under the translations, which is also the conditional e-variable of \citet[Remark~6]{hao2026evalues}, attains this value~\citep[Corollary~5]{perezortiz2024estatistics}.
These likelihood ratios form a test martingale with respect to the filtration of the maximal invariants~\citep[Proposition~6]{perezortiz2024estatistics}.
However, they do not form an e-process with respect to the data filtration (Proposition~\ref{prop:notep}).
Indeed, against the null distribution $\N(\theta,1)^{\otimes\infty}$, a stopping time may depend on the distance between the running mean and $\theta$, which is not determined by the maximal invariants.
Our counterexample is thus a specific case of group-invariant testing, in which the law of the iterated logarithm for the running mean determines the regret beyond $\kappa$.
We prove in Theorem~\ref{thm:composite} that, over all horizons up to $T$, the minimax regret is $\kappa+\{1-s^2+o(1)\}\log\log T$, and that a single explicit e-process attains this rate at every horizon.
In contrast, if $s\ge1$, the minimax regret is at most $\log s$ for every $T$ (Proposition~\ref{prop:sge1} and Section~\ref{sec:rate}).
Thus, the minimax regret is bounded if and only if the alternatives have at least the variance of the null.
The growth of the minimax regret provides a counterexample to general bounded-regret attainability for the question of \citet[Section~5]{grunwald2024safe}.
It also refutes the conjecture of \citet{hao2026evalues} if ``almost'' means that their conditional e-variables have bounded expectations at stopping times.
In contrast with the statement of \citet{koning2026anytime}, which is concerned with tests, the anytime validity of e-processes is therefore not free when the optimal e-power is required at every horizon up to $T$.

\section{Setting and basic bounds}\label{sec:setting}

The numéraire $\Es_t$ of~\eqref{eq:numeraire} is positive and unique up to $Q$-null sets \citep[Definition~2.1 and Proposition~2.4]{larsson2025numeraire}, and $\Es_0=1$ since $\cF_0$ is trivial.
The measure $\Ps_t$ is the reverse information projection of $Q_t$ on the convex hull of the null distributions of $X^t$ whenever this projection exists \citetext{\citealp{grunwald2024safe}; \citealp[Corollaries~4.10 and~4.11]{larsson2025numeraire}}.
Because the second identity in~\eqref{eq:numeraire} requires a finite right-hand side, we assume throughout that $\sup_{E\in\cE_t}\E_Q[\log E]<\infty$ for every $t\ge0$.
It follows from Lemma~\ref{lem:regret} that $\log\Es_t$ is $Q$-integrable and that $R_t(E)=\E_Q[\log\Es_t]-\E_Q[\log E]\in[0,\infty]$ for every $E\in\cE_t$.
For $T\ge0$, we fix versions of $\Es_0,\ldots,\Es_T$ and let
\begin{equation}\label{eq:snell}
  \eta_T=\sup_{P\in\cP}\sup_\tau\E_P[\Es_{\tau\wedge T}],
\end{equation}
where $\tau$ ranges over the stopping times of $(\cF_t)$ and $a\wedge b=\min(a,b)$.
For $P\in\cP$, the Snell envelope~\citep[Section~3.3]{ramdas2022testing} of $(\Es_{t\wedge T})_{t\ge0}$ under $P$ is the process $(\Lambda^P_t)_{0\le t\le T}$ defined by
\[
  \Lambda^P_T=\Es_T,\quad\Lambda^P_t=\max\{\Es_t,\E_P[\Lambda^P_{t+1}\mid\cF_t]\}
\]
for $t<T$.
It is the smallest supermartingale under $P$ with respect to $(\cF_t)$ that dominates $(\Es_{t\wedge T})_{0\le t\le T}$~\citep[Theorem~1.2]{peskir2006optimal}.
Since $\cF_0$ is trivial, the same theorem gives $\Lambda^P_0=\sup_\tau\E_P[\Es_{\tau\wedge T}]$, and hence $\eta_T=\sup_{P\in\cP}\Lambda^P_0$.

Let $\Delta_T$ be the set of probability vectors $(\xi_0,\ldots,\xi_T)$.
To combine the numéraires at the horizons $0,\ldots,T$ into one e-process, we weight them by a vector $\xi\in\Delta_T$.
The following proposition gives the resulting bounds on the regret.
\begin{proposition}[Basic bounds]\label{prop:basic}
  Let $T\ge0$.
  Then $\Gamma_T\le\Gamma_{T+1}$.
  \begin{parts}
  \item\partlabel{prop:basic:i} For every $t\ge0$, the process $(E_i)_{i\ge0}$ with $E_i=\Es_t\ind(i\ge t)$ is an e-process with $R_t(E_t)=0$, where $\ind(\cdot)$ is the indicator function.
  \item\partlabel{prop:basic:ii} For every $\xi\in\Delta_T$, the process $(\sum_{i=0}^{t\wedge T}\xi_i\Es_i)_{t\ge0}$ is an e-process with regret at most $\log(1/\xi_t)$ at every $t\le T$.
    Consequently, $\Gamma_T\le\log(T+1)$.
  \item\partlabel{prop:basic:iii} If $\eta_T<\infty$, then the process $(E_t)_{t\ge0}$ with $E_t=\Es_t\ind(t\le T)/\eta_T$ is an e-process with regret $\log\eta_T$ at every $t\le T$.
    Hence, $\Gamma_T\le\log\eta_T$.
  \end{parts}
\end{proposition}

\begin{remark}\label{rem:zero}
  If $\Es_0,\ldots,\Es_T$ have versions such that $(\Es_{t\wedge T})_{t\ge0}$ is an e-process, then $\Gamma_T=0$.
  Such versions exist for a simple null~\citep[Theorem~8]{koolen2022log} and for some exponential families with independent and identically distributed data (Section~\ref{sec:finite}).
  Conversely, $\Gamma_T=0$ only if such versions exist, when $\cX$ is finite or the null distributions of $X^T$ are absolutely continuous with respect to $Q_T$ (Proposition~\ref{prop:minimax:i}).
  In the Bernoulli example of \citet[Section~4.5]{koolen2022log}, the numéraires do not form an e-process, but $\Gamma_T$ appears empirically bounded (Example~\ref{ex:bernoulli}).
\end{remark}

\section{A lower bound by duality}\label{sec:duality}

We now combine the primal-dual equivalence of \citet[Theorem~5]{koolen2022log} with the optimal e-powers at fixed horizons to obtain a lower bound on $\Gamma_T$.
One of its two inequalities bounds the average over the horizons of the e-powers of an e-process by the Kullback--Leibler divergence of the stopped alternative from a mixture of null distributions stopped at random times.

The stopping times in these mixtures may be randomized.
A randomized stopping time is defined as a measurable map $\tau\colon\cX^\infty\times[0,1]\to\{0,1,\ldots\}\cup\{\infty\}$ such that $\tau(\cdot,u)$ is a stopping time of $(\cF_t)$ for every $u\in[0,1]$.
Applied to each $\tau(\cdot,u)$, the second inequality in~\eqref{eq:eproc} gives $\int_0^1\E_P[E_{\tau(X,u)}]\dd u\le1$ for every e-process $(E_t)$, every $P\in\cP$, and every randomized stopping time $\tau$.

Distributions stopped at a random time are defined on the pairs $(t,x^t)$ with $t\in\{0,1,\ldots\}\cup\{\infty\}$ and $x^t\in\cX^t$.
For $\xi\in\Delta_T$, let $Q^\xi[\{t\}\times B]=\xi_tQ_t(B)$ for $0\le t\le T$.
In other words, $Q^\xi$ is the distribution of the path of $X\sim Q$ stopped at an independent time with distribution $\xi$.
For $P\in\cP$ and a stopping time $\tau$ of $(\cF_t)$, let $P^\tau$ be the distribution of $(\tau,X^\tau)$ when $X\sim P$, with $X^\tau=(X_t)_{t\ge1}$ on the event $\{\tau=\infty\}$.
For a randomized stopping time $\tau$, let $P^\tau=\int_0^1P^{\tau(\cdot,u)}\dd u$.
The set $\cB$ of mixtures of stopped null distributions consists of the distributions
\[
  \bar P=\int P_\lambda^{\tau_\lambda}\pi(\dd\lambda),
\]
where $\pi$ is a probability measure on a measurable space, $\lambda\mapsto P_\lambda\in\cP$ is a probability kernel, and $(\lambda,x,u)\mapsto\tau_\lambda(x,u)$ is measurable with $\tau_\lambda$ a randomized stopping time for every $\lambda$.
The set $\cB$ contains the Bayesian mixtures of stopped elements of the null of \citet[Section~3]{koolen2022log}.
For finite $\cX$, \citet[Theorem~5]{koolen2022log} show that the supremum of $\sum_t\xi_t\E_Q[\log E_t]$ over e-processes $(E_t)$ equals the infimum of $\KL(Q^\xi\|\bar P)$ over their Bayesian mixtures $\bar P$, where $\KL$ is the Kullback--Leibler divergence.
More generally, the inequality $\sum_t\xi_t\E_Q[\log E_t]\le\KL(Q^\xi\|\bar P)$ holds for every $\cX$ and every $\bar P\in\cB$, as we prove in Section~\ref{sec:proofs}.
Together with Lemma~\ref{lem:regret}, it gives the following lower bound.

\begin{lemma}[Lower bound by duality]\label{lem:duality}
  Let $T\ge0$.
  For every e-process $(E_t)$, $\xi\in\Delta_T$, and $\bar P\in\cB$,
  \begin{equation}\label{eq:duality}
    \sum_{t=0}^T\xi_tR_t(E_t)\ge\sum_{t=0}^T\xi_t\E_Q[\log\Es_t]-\KL(Q^\xi\|\bar P),
  \end{equation}
  with the convention $0\times\infty=0$ in the sums.
  Consequently, for every $\xi\in\Delta_T$,
  \begin{equation}\label{eq:gamma-lb}
    \Gamma_T\ge\inf_{(E_t)}\sum_{t=0}^T\xi_tR_t(E_t)\ge\sum_{t=0}^T\xi_t\E_Q[\log\Es_t]-\inf_{\bar P\in\cB}\KL(Q^\xi\|\bar P).
  \end{equation}
\end{lemma}

\begin{remark}\label{rem:deterministic}
  If every $\tau_\lambda(x,u)$ depends on $u$ only, then $\KL(Q^\xi\|\bar P)\ge\sum_t\xi_t\E_Q[\log\Es_t]$, and the right-hand side of~\eqref{eq:duality} is at most zero.
\end{remark}

\section{The Gaussian location model}\label{sec:gauss}

\subsection{Numéraires}\label{sec:numeraires}

Throughout this section, $\cX=\R$, $s>0$, $\rho\ge0$, and $\rho^2\ge1-s^2$.
The null and the alternative are
\[
  \cP=\{P_\theta:\theta\in\R\},\quad P_\theta=\N(\theta,1)^{\otimes\infty},\quad Q=\int Q_\mu\N(0,\rho^2)(\dd\mu),\quad Q_\mu=\N(\mu,s^2)^{\otimes\infty},
\]
where $\N(0,0)$ is the unit mass at zero by convention.
The alternative of interest is the composite alternative $\{Q_\mu:\mu\in\R\}$.
For this alternative, we take the regret in the worst case over $\mu$.
This worst case is at least the average of the regret under $\N(0,\rho^2)$.
We therefore study the mixture $Q$ first and return to the worst case over $\mu$ in Theorem~\ref{thm:composite}.
We write $\E_\theta$ for the expectation under $P_\theta$, and $P_{\theta,t}$ and $Q_{\mu,t}$ for the distributions of $X^t$ under $P_\theta$ and $Q_\mu$.
Let $\bar X_t=t^{-1}\sum_{i=1}^tX_i$ be the running mean, let $\varphi_v$ be the density of $\N(0,v)$, and define $Y_2,Y_3,\ldots$ by
\begin{equation}\label{eq:Y}
  Y_i=\{(i-1)/i\}^{1/2}(X_i-\bar X_{i-1}).
\end{equation}
We write $Y^t=(Y_2,\ldots,Y_t)$, and we let $(\cG_t)_{t\ge0}$ be the filtration generated by $(Y_t)_{t\ge2}$, with $\cG_0$ and $\cG_1$ trivial.
The linear map from $X^t$ to $(t^{1/2}\bar X_t,Y_2,\ldots,Y_t)$ is orthogonal~\citep{lancaster1965helmert}, that is, $\sum_{i=1}^t(X_i-m)^2=t(\bar X_t-m)^2+\sum_{i=2}^tY_i^2$ for every $m\in\R$ (Lemma~\ref{lem:orth}).
It follows that the variables $\bar X_t,Y_2,\ldots,Y_t$ are independent under every $P_\theta$ and every $Q_\mu$, with $Y_i\sim\N(0,1)$ under $P_\theta$ and $Y_i\sim\N(0,s^2)$ under $Q_\mu$.
The null and alternative thus differ in the variance of the $Y_i$, while only $\bar X_t$ depends on the location.
We show that $\Ps_t$ of~\eqref{eq:numeraire} is a mixture of null distributions under which $\bar X_t$ has the same distribution as under $Q_t$.

\begin{proposition}[Numéraires of the Gaussian location model]\label{prop:numeraire}
  For every $t\ge1$,
  \[
    \Es_t=\prod_{i=2}^t\ell(Y_i),\quad\Ps_t=\int P_{\theta,t}\N\{0,\rho^2-(1-s^2)/t\}(\dd\theta),\quad\E_Q[\log\Es_t]=(t-1)\kappa,
  \]
  where $\kappa=\KL\{\N(0,s^2)\|\N(0,1)\}=(s^2-1-\log s^2)/2$ and $\ell=\varphi_{s^2}/\varphi_1$.
  Under every $P_\theta$, $(\Es_t)_{t\ge0}$ is a test martingale with respect to $(\cG_t)_{t\ge0}$.
\end{proposition}

The condition $\rho^2\ge1-s^2$ makes the mixing variance $\rho^2-(1-s^2)/t$ of $\Ps_t$ nonnegative for every $t\ge1$.
Since the vector $Y^t$ is a maximal invariant under translations $x\mapsto(x_1+m,\ldots,x_t+m)$, the numéraire $\Es_t$ is the likelihood ratio of a maximal invariant~\citep{perezortiz2024estatistics}.
It is also the conditional e-variable of \citet{hao2026evalues}.
We note that the formulas for $\Ps_t$ and $\E_Q[\log\Es_t]$ in Proposition~\ref{prop:numeraire} are the one-dimensional case of their Theorems~1 and~2.

Under $P_\theta$ and given $X_1$, the variable $Y_2=(X_2-X_1)/2^{1/2}$ is normal with mean $(\theta-X_1)/2^{1/2}$ and variance $1/2$ instead of $1$.
For $s<1$, $\ell(y)$ decreases in $|y|$, and its expectation under $\N(0,1)$ is one.
The conditional expectation of $\Es_2=\ell(Y_2)$ therefore exceeds one when $X_1$ is close to $\theta$.
For $s<1$, the numéraires are therefore not an e-process with respect to $(\cF_t)$.
Proposition~\ref{prop:notep} gives the conditional expectation of $\Es_2$ given $X_1$ and the best stopping time up to $t=2$.

\subsection{Stopping rules that depend on the null mean}\label{sec:stopping}

Since validity is required under each null distribution separately, a different stopping time may be used against each $P_\theta$.
For $\theta\in\R$, let $L^\theta_t=\prod_{i=1}^t\ell(X_i-\theta)$ be the likelihood ratio of $Q_\theta$ against $P_\theta$ on $\cF_t$, with $L^\theta_0=1$.
Under $P_\theta$, $(L^\theta_t)_{t\ge0}$ is a nonnegative martingale with respect to $(\cF_t)$ with initial value one.

\begin{lemma}[Change of measure]\label{lem:identity}
  \begin{parts}
  \item\partlabel{lem:identity:i} For every $\theta\in\R$ and every $t\ge1$, $\Es_t=sL^\theta_t\exp\{(1-s^2)t(\bar X_t-\theta)^2/(2s^2)\}$.
  \item\partlabel{lem:identity:ii} Let $(S_t)_{t\ge0}$ be a random walk with $S_0=0$ and independent $\N(0,1)$ increments, and let $h_0=1$ and $h_t=s\exp\{(1-s^2)S_t^2/(2t)\}$ for $t\ge1$.
    For every $T\ge0$, $\eta_T=\sup_\tau\E[h_{\tau\wedge T}]$, where $\tau$ ranges over the stopping times of the filtration generated by $(S_t)$.
  \end{parts}
\end{lemma}

Lemma~\ref{lem:identity:i} splits the numéraire into the likelihood ratio $L^\theta_t$ and a function of the distance between the running mean and $\theta$.
Since $L^\theta_t$ is a martingale under $P_\theta$, one can raise $\E_\theta[\Es_\tau]$ above one only through the second factor.
Part~\ref{lem:identity:ii*} then turns $\eta_T$ into an optimal stopping problem for a Gaussian random walk.
If $s<1$, the reward $h_t$ is large when $|S_t|$ is large compared with $t^{1/2}$.
By the law of the iterated logarithm, $\limsup_{t\to\infty}\log h_t/\log\log t=1-s^2$ almost surely.
In Section~\ref{sec:rate}, we show that this limit determines the rate $(1-s^2)\log\log T$ of both $\log\eta_T$ and $\Gamma_T$.

\subsection{The iterated-logarithm rate}\label{sec:rate}

The following proposition gives the rate of $\Gamma_T$ for the mixture and $s<1$, and an e-process that attains it at every horizon.
For the lower bound, we apply Lemma~\ref{lem:duality} to a mixture over $\theta$ of null distributions stopped at the first horizon of a geometric grid at which $t^{1/2}|\bar X_t-\theta|$ exceeds the level $s(2\log\log T)^{1/2}$.
To bound the density ratio of Lemma~\ref{lem:duality} from below, half of the mixture stops instead at an independent uniform horizon of the grid.
At each horizon of the grid, $t^{1/2}|\bar X_t-\theta|$ exceeds $s(2\log\log T)^{1/2}$ with a probability of order $1/\log T$ under $Q_\theta$ and of order $(\log T)^{-s^2}$ under $P_\theta$.
The logarithm of the ratio of these probabilities gives the lower bound $(1-s^2)\log\log T$.
For the upper bound, we control the maximum of $(h_t)$ on geometric blocks of times by Doob's maximal inequality (Lemma~\ref{lem:upper}).
\begin{proposition}[Iterated-logarithm rate for the mixture]\label{prop:mixture}
  Let $s\in(0,1)$ and $\rho^2\ge1-s^2$.
  \begin{parts}
  \item\partlabel{prop:mixture:i} As $T\to\infty$, $\Gamma_T/\log\log T\to1-s^2$ and $\log\eta_T/\log\log T\to1-s^2$.
  \item\partlabel{prop:mixture:ii} For $j\ge1$, let $T_j=\lceil\exp(e^j)\rceil$ and $q_j=1+s^2/\{(1-s^2)(j+1)\}$.
    The process $(E_t)_{t\ge0}$ with
    \[
      E_t=\Es_t\sum_{j:T_j\ge t}\frac{1}{j(j+1)u_j},\quad u_j=\frac{j+1}{s^2j}\biggl(2+\frac{2\log T_j}{\log q_j}\biggr)^{1-s^2j/(j+1)},
    \]
    is an e-process with $\limsup_{t\to\infty}R_t(E_t)/\log\log t=1-s^2$.
  \end{parts}
\end{proposition}
\begin{figure}
  \centering
\includegraphics{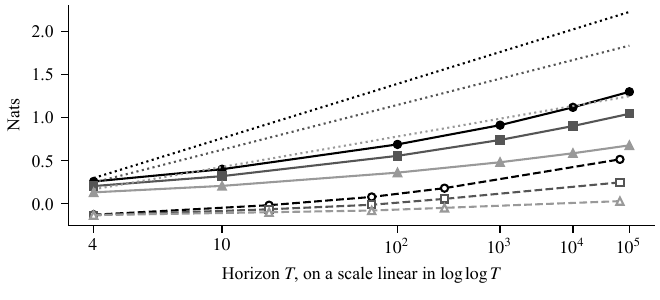}
  \caption{Gaussian location model with $s=0.3$, $0.5$, and $0.7$: upper bound on $\log\eta_T$, computed by backward induction (\legEtaThree, \legEtaFive, \legEtaSeven); Monte Carlo estimate of the limit as $\rho\to\infty$ of the lower bound in~\eqref{eq:S-lb} at the horizons $T=4^K$ (\legLowThree, \legLowFive, \legLowSeven); and the function $(1-s^2)\log\log T$ of Proposition~\ref{prop:mixture:i} (\legLimThree, \legLimFive, \legLimSeven). The first two bound the minimax regret in the worst case over $\mu$ minus $\kappa$ from above and from below.}
  \label{fig:bounds}
\end{figure}

Proposition~\ref{prop:mixture:i} holds for every $\rho>0$ (Proposition~\ref{prop:smallrho}).
For $s<1$ and a single horizon $t\ge1$, the numéraire $\Es_t$ attains the minimax regret $\kappa$ in the worst case over $\mu$ (see the proof of Theorem~\ref{thm:composite}).
The following theorem, our main result, shows that validity at all horizons up to $T$ adds between $\Gamma_T$ and $\log\eta_T$ to $\kappa$.
In contrast, for each single $Q_\mu$, the minimax regret over all horizons is zero (Proposition~\ref{prop:simple:i}).
\begin{theorem}[Composite alternative]\label{thm:composite}
  Let $s\in(0,1)$ and $\rho^2\ge1-s^2$.
  For $\mu\in\R$, let $R^\mu_t$ be the regret defined by~\eqref{eq:numeraire} and~\eqref{eq:price} with $Q_\mu$ in place of $Q$.
  \begin{parts}
  \item\partlabel{thm:composite:i} For $T\ge1$,
    \[
      \kappa+\Gamma_T\le\inf_{(E_t)}\max_{0\le t\le T}\sup_{\mu\in\R}R^\mu_t(E_t)\le\kappa+\log\eta_T.
    \]
    The middle term divided by $\log\log T$ therefore tends to $1-s^2$ as $T\to\infty$.
  \item\partlabel{thm:composite:ii} The e-process $(E_t)$ of Proposition~\ref{prop:mixture:ii} satisfies $\sup_{\mu\in\R}R^\mu_t(E_t)=\kappa+R_t(E_t)$ for $t\ge1$, so that $\limsup_{t\to\infty}\sup_{\mu\in\R}R^\mu_t(E_t)/\log\log t=1-s^2$.
  \end{parts}
\end{theorem}
Theorem~\ref{thm:composite} shows that the minimax regret in the worst case over $\mu$ is unbounded.
It grows more slowly than any power of $\log T$, while the optimal e-power $T\kappa$ of each $Q_\mu$ grows linearly.
The terms of the e-process in~\ref{thm:composite:ii*} are the numéraires multiplied by a deterministic factor that decreases to zero.
This e-process therefore attains the rate of~\ref{thm:composite:i*} without a planned horizon.
Since $\ell(y)=s^{-1}\exp\{-(1-s^2)y^2/(2s^2)\}$, the numéraire $\Es_t=\prod_{i=2}^t\ell(Y_i)$ is a function of $\sum_{i=2}^tY_i^2=\sum_{i=1}^t(X_i-\bar X_t)^2$ (Lemma~\ref{lem:orth:i}).
The terms of the e-process in~\ref{thm:composite:ii*} are therefore invariant under permutations of $X_1,\ldots,X_t$.
As shown in Fig.~\ref{fig:bounds}, at $T=10^5$, $\log\eta_T$ is still below three fifths of $(1-s^2)\log\log T$ for each of the three values of $s$.
The convergence in Theorem~\ref{thm:composite:i} is therefore slow.
Finally, the following proposition treats the case $s\ge1$.
\begin{proposition}\label{prop:sge1}
  If $s\ge1$, then $\Gamma_T\le\log\eta_T\le\log s$ for every $T\ge0$, and $\log\eta_T\to\log s$ as $T\to\infty$.
\end{proposition}
For $s\ge1$, the regret is also bounded in the worst case over $\mu$.
Indeed, $\prod_{i=2}^t\ell(Y_i)$ of Proposition~\ref{prop:numeraire} is then also the numéraire of every $Q_\mu$ (Proposition~\ref{prop:simple}).
The e-process of Proposition~\ref{prop:basic:iii} therefore has regret at most $\log s$ against every $Q_\mu$ at every $t\le T$.
Nevertheless, $\Gamma_T>0$ for $s>1$ and $T\ge2$ (Section~\ref{sec:notep}).

\section{Discussion}\label{sec:discussion}

From a theoretical point of view, the unbounded minimax regret is disappointing, since no e-process attains the optimal e-power at every horizon up to a bounded loss.
For practitioners, however, the size of this loss is hardly a concern.
In the Gaussian location model, the optimal e-power grows by $\kappa=(s^2-1-\log s^2)/2$ per observation, where $s^2$ is the variance of the alternatives and one the variance of the null (Proposition~\ref{prop:numeraire}).
A regret of $\log\eta_T$ thus equals the optimal e-power of $\log\eta_T/\kappa$ observations.
This equality holds for e-power only, not for the expected stopping time or the sample size needed for a given power.
By Proposition~\ref{prop:mixture:i}, the number $\log\eta_T/\kappa$ of observations is asymptotically $(1-s^2)\log\log T/\kappa$, and $\log\log T<4$ for $T<10^{23}$.
For $s$ close to one, $\kappa\approx(1-s^2)^2/4$, so that this number grows like $4\log\log T/(1-s^2)$.
The e-process with regret $\log\eta_T$ at every horizon up to $T$ is the conditional e-variable of \citet{hao2026evalues} divided by $\eta_T$ up to $T$ and equal to zero after $T$ (Proposition~\ref{prop:basic:iii}).
An upper bound on $\eta_T$, from Table~\ref{tab:snell} or Lemma~\ref{lem:upper}, may replace $\eta_T$.
The division by $\eta_T$ is needed only when the value at a data-dependent time must remain an e-variable, for instance to be combined with later evidence.
In contrast, the test that rejects the null as soon as the conditional e-variables reach $1/\alpha$ has level $\alpha$, while the e-process reaches $1/\alpha$ only when they reach $\eta_T/\alpha$.
Indeed, these e-variables form a test martingale with respect to the filtration $(\cG_t)$ generated by the variables $Y_i$ of~\eqref{eq:Y}.
Therefore, only the e-value property at data-dependent times makes the minimax regret grow.

Other general constructions of e-processes have a bounded regret only under stronger assumptions, or lose more e-power.
The universal-inference e-process~\citep[Section~7]{wasserman2020universal} has a bounded regret for an exponential-family null and a fixed independent and identically distributed alternative~\citep[Theorem~3 and Corollary~3]{hao2026evalues}.
Replacing the null mean by the running mean in the likelihood ratios~\citep[Theorem~2, Part~3]{hao2026evalues} gives a regret of order $\log t$ in the worst case over the location for $s<1$ (Proposition~\ref{prop:plugin}).
The adjuster of \citet[Theorem~2]{choe2026combining}, applied to the numéraires, gives a simulated regret that grows like $2\log t$ (Section~\ref{sec:numerics}).

Iterated logarithms also appear in lower bounds on expected stopping times \citetext{\citealp[Theorem~1]{farrell1964asymptotic}; \citealp[Corollary~3.3]{agrawal2025stopping}} and in the estimation of a Gaussian mean~\citep[Theorem~3.2]{kirichenko2020minimax}.
In our model, the stopping times of the lower bound depend on the running mean, the fluctuations of which are governed by the law of the iterated logarithm.
We therefore expect a minimax regret of order $\log\log T$ in a larger class of models.
Such models may include all group-invariant problems in which the likelihood ratio of a maximal invariant attains the minimax regret at each single horizon, such as the $t$-test \citetext{\citealp[Corollary~5 and Example~1]{perezortiz2024estatistics}; \citealp{wang2025anytime,grunwald2026supermartingales}}.
Two further questions remain open.
First, for $s>1$, it is unknown whether $\Gamma_T$ tends to $\log s$, the limit of $\log\eta_T$.
Second, the pairs $(\cP,Q)$ for which $\Gamma_T\to\infty$ are not characterized.

\section*{Use of generative AI}

During the preparation of this work the author used Claude (Anthropic) in order to revise the text and check proofs, and to write and run the numerical code.
The author also used Claude to formalize all results in the Lean proof assistant (see Section~\ref{sec:formalization} for details).
As a result, Proposition~\ref{prop:minimax:ii} was strengthened, Proposition~\ref{prop:smallrho} was extended to the limit of $\log\eta_T$, and the proofs of Lemma~\ref{lem:duality}, Remark~\ref{rem:deterministic}, and Proposition~\ref{prop:minimax:i} were shortened.

\section*{Acknowledgement}

Unbeknownst to the author, a claim that the regret under the same condition $s>1$ is bounded was given by Stephan Bongers and Peter Gr\"unwald in the poster ``When generalized likelihood ratios are---and are not---e-processes'' at the SAVI Conference in July 2026 (\url{https://savi-conference.github.io/2026/program/}).
The poster did not explicitly give a proof, but showed simulation results.
The growth of the regret when $s < 1$ did not appear there.
The author thanks Peter Gr\"unwald for pointing this out and for a question that led to the present form of the e-process of Proposition~\ref{prop:mixture:ii} (see Remark~\ref{rem:stopped}).

\bibliographystyle{plainnat}
\bibliography{refs}

\clearpage
\appendix
\setcounter{section}{0}\def\thesection{S\arabic{section}}\def\theHsection{S\arabic{section}}
\setcounter{equation}{0}\def\theequation{S\arabic{equation}}\def\theHequation{S\arabic{equation}}
\setcounter{table}{0}\def\thetable{S\arabic{table}}\def\theHtable{S\arabic{table}}
\setcounter{figure}{0}\def\thefigure{S\arabic{figure}}\def\theHfigure{S\arabic{figure}}
\setcounter{theorem}{0}\def\thetheorem{S\arabic{theorem}}\def\theHtheorem{S\arabic{theorem}}
\setcounter{proposition}{0}\def\theproposition{S\arabic{proposition}}\def\theHproposition{S\arabic{proposition}}
\setcounter{lemma}{0}\def\thelemma{S\arabic{lemma}}\def\theHlemma{S\arabic{lemma}}
\setcounter{corollary}{0}\def\thecorollary{S\arabic{corollary}}\def\theHcorollary{S\arabic{corollary}}
\setcounter{definition}{0}\def\thedefinition{S\arabic{definition}}\def\theHdefinition{S\arabic{definition}}
\setcounter{example}{0}\def\theexample{S\arabic{example}}\def\theHexample{S\arabic{example}}
\setcounter{remark}{0}\def\theremark{S\arabic{remark}}\def\theHremark{S\arabic{remark}}

\begin{center}
\Large Supplementary material
\end{center}

\section{Proofs}\label{sec:proofs}

Throughout, expectations are defined in $[-\infty,\infty]$ whenever the positive part or the negative part is integrable.
We define $\log^+x=\max(\log x,0)$ and $\log^-x=\max(-\log x,0)$ to be the positive and negative parts of $\log x$ for $x\in[0,\infty]$.
We will use the fact that $\log^-x\le1/x$ several times.

\subsection{Technical results}\label{sec:technical}

The following lemma gives the regret as the loss of e-power relative to the numéraire.

\begin{lemma}[Integrability of the numéraire]\label{lem:regret}
  For every $t\ge0$, $\log\Es_t$ is $Q$-integrable and satisfies the second identity in~\eqref{eq:numeraire}.
  For every $E\in\cE_t$, the regret $R_t(E)$ in~\eqref{eq:price} is well defined, and $R_t(E)=\E_Q[\log\Es_t]-\E_Q[\log E]\in[0,\infty]$.
\end{lemma}

\begin{proof}[of Lemma~\ref{lem:regret}]
  Taking $E\equiv1$ in the first property in~\eqref{eq:numeraire} gives $\E_Q[1/\Es_t]\le1$.
  Since $\log^-x\le1/x$, it follows that $\E_Q[\log^-\Es_t]\le1$.
  Hence, $\E_Q[\log\Es_t]$ is well defined, and
  \[
    \E_Q[\log\Es_t]\le\sup_{E\in\cE_t}\E_Q[\log E]<\infty.
  \]
  Therefore, $\log\Es_t$ is integrable, and $0<\Es_t<\infty$ $Q$-almost surely.

  Let $E\in\cE_t$.
  Since $\log^+(E/\Es_t)\le E/\Es_t$, the expectation $\E_Q[\log(E/\Es_t)]$ is well defined in $[-\infty,1]$.
  Jensen's inequality gives
  \[
    \E_Q[\log(E/\Es_t)]\le\log\E_Q[E/\Es_t]\le0.
  \]
  Since $\log E=\log\Es_t+\log(E/\Es_t)$ $Q$-almost surely, with an integrable first term,
  \[
    \E_Q[\log E]=\E_Q[\log\Es_t]-R_t(E)\le\E_Q[\log\Es_t].
  \]
  Taking the supremum over $E\in\cE_t$ gives the second identity in~\eqref{eq:numeraire}, and the formula for $R_t(E)$ follows.
\end{proof}

The following lemma gives standard facts on the variables $Y_i$ of~\eqref{eq:Y}~\citep{lancaster1965helmert}.

\begin{lemma}[Orthogonal transformation]\label{lem:orth}
  Let $t\ge2$, $x\in\R^t$, $\bar x_i=i^{-1}\sum_{j=1}^ix_j$ for $1\le i\le t$, and $y_i=\{(i-1)/i\}^{1/2}(x_i-\bar x_{i-1})$ for $2\le i\le t$.
  \begin{parts}
  \item\partlabel{lem:orth:i} The map $x\mapsto(t^{1/2}\bar x_t,y_2,\ldots,y_t)$ is orthogonal.
    Moreover, the values $y_2,\ldots,y_t$ do not change when the same constant is added to every $x_i$.
    Hence, $\sum_{i=1}^t(x_i-m)^2=t(\bar x_t-m)^2+\sum_{i=2}^ty_i^2$ for every $m\in\R$.
  \item\partlabel{lem:orth:ii} $\bar x_t=\bar x_{t-1}+y_t/\{t(t-1)\}^{1/2}$.
  \item\partlabel{lem:orth:iii} If $X^t\sim\N(m,v)^{\otimes t}$, then $\bar X_t,Y_2,\ldots,Y_t$ are independent, with $\bar X_t\sim\N(m,v/t)$ and $Y_i\sim\N(0,v)$.
  \end{parts}
\end{lemma}

\begin{proof}[of Lemma~\ref{lem:orth}]
  Let $\omega_1=t^{-1/2}(1,\ldots,1)$ and, for $2\le i\le t$, let
  \[
    \omega_i=\{i(i-1)\}^{-1/2}(1,\ldots,1,1-i,0,\ldots,0),
  \]
  with $i-1$ leading ones.
  These vectors are the rows of a Helmert matrix~\citep{lancaster1965helmert} and form an orthonormal basis of $\R^t$, with
  \[
    \langle\omega_1,x\rangle=t^{1/2}\bar x_t,\quad\langle\omega_i,x\rangle=-y_i.
  \]
  Therefore, the map in~\ref{lem:orth:i*} is orthogonal.
  Since $\omega_i$ is orthogonal to $(1,\ldots,1)$ for $i\ge2$, the values $y_2,\ldots,y_t$ do not change when a constant is added to every coordinate.
  Applying the orthogonality to $(x_1-m,\ldots,x_t-m)$ gives the identity for $\sum_i(x_i-m)^2$ in part~\ref{lem:orth:i*}.
  The identity $\bar x_t=\bar x_{t-1}+(x_t-\bar x_{t-1})/t$ gives part~\ref{lem:orth:ii*}.
  Since the image of a normal distribution with covariance $v$ times the identity under an orthogonal map is normal with the same covariance, part~\ref{lem:orth:i*} gives part~\ref{lem:orth:iii*}.
\end{proof}

In the model of Section~\ref{sec:gauss}, with $b=s^{-2}-1$, the likelihood ratio $\ell$ satisfies $\log\ell(y)=-\log s-by^2/2$.
For $m\in\R$, $v\ge0$ with $1+bv>0$, and $Z\sim\N(0,1)$, completing the square gives
\begin{equation}\label{eq:S-gauss}
  \E[\ell(m+v^{1/2}Z)]=\{s(1+bv)^{1/2}\}^{-1}\exp\{-bm^2/(2+2bv)\}.
\end{equation}
Since $s(1+b)^{1/2}=1$ and $b/(1+b)=1-s^2$, the right-hand side of~\eqref{eq:S-gauss} equals $\exp\{-(1-s^2)m^2/2\}$ for $v=1$.
By Lemma~\ref{lem:orth:i}, for $x\in\R^t$, every $m\in\R$, and every $v>0$,
\begin{equation}\label{eq:S-factor}
  \prod_{i=1}^t\varphi_v(x_i-m)=t^{-1/2}\varphi_{v/t}(\bar x_t-m)\prod_{i=2}^t\varphi_v(y_i).
\end{equation}

\subsection{Proof of Proposition~\ref{prop:basic}}\label{sec:proof-basic}

\begin{proof}[of Proposition~\ref{prop:basic}]
  Since the maximum in~\eqref{eq:price} is taken over more horizons for $T+1$ than for $T$, $\Gamma_T\le\Gamma_{T+1}$.

  (i) Fix $t\ge0$ and let $E_i=\Es_t\ind(i\ge t)$ for $i\ge0$.
  For every $P\in\cP$ and every stopping time $\tau$ of $(\cF_t)$,
  \[
    \E_P[E_\tau]=\E_P[\Es_t\ind(\tau\ge t)]\le\E_P[\Es_t]\le1.
  \]
  Hence, $(E_i)_{i\ge0}$ is an e-process with $R_t(E_t)=0$.

  (ii) Let $E_t=\sum_{i=0}^{t\wedge T}\xi_i\Es_i$.
  Since $(E_t)_{t\ge0}$ is stopped at $T$, every $P\in\cP$ and every stopping time $\tau$ of $(\cF_t)$, including the case $\tau=\infty$, satisfy
  \[
    \E_P[E_\tau]=\sum_{i=0}^T\xi_i\E_P[\Es_i\ind(\tau\ge i)]\le\sum_{i=0}^T\xi_i=1.
  \]
  Since $E_t\ge\xi_t\Es_t$, the regret satisfies $R_t(E_t)\le\log(1/\xi_t)$.
  In particular, the uniform $\xi$ gives $\Gamma_T\le\log(T+1)$.

  (iii) Since $\Es_0=1$, taking $\tau=0$ in~\eqref{eq:snell} gives $\eta_T\ge1$.
  If $\eta_T=\infty$, the bound $\Gamma_T\le\log\eta_T$ holds trivially.
  Let $\eta_T<\infty$, and let $E_t=\Es_t\ind(t\le T)/\eta_T$ for $t\ge0$, with the versions of $\Es_0,\ldots,\Es_T$ fixed in~\eqref{eq:snell}.
  Then $E_t$ is $\cF_t$-measurable, with values in $[0,\infty]$.
  Since $E_t=0$ for every $t>T$, we have $E_\infty=0$.
  Thus, $E_\tau=\Es_{\tau\wedge T}\ind(\tau\le T)/\eta_T$ for every stopping time $\tau$ of $(\cF_t)$, also on the event $\{\tau=\infty\}$.
  It follows from~\eqref{eq:snell} that every $P\in\cP$ and every such $\tau$ satisfy
  \[
    \E_P[E_\tau]\le\E_P[\Es_{\tau\wedge T}]/\eta_T\le1.
  \]
  Consequently, $(E_t)_{t\ge0}$ is an e-process.

  Let $t\le T$.
  Since $\log\Es_t$ is $Q$-integrable by Lemma~\ref{lem:regret}, $0<\Es_t<\infty$ $Q$-almost surely.
  Therefore, $\Es_t/E_t=\eta_T$ $Q$-almost surely, and
  \[
    R_t(E_t)=\E_Q[\log(\Es_t/E_t)]=\log\eta_T.
  \]
  Taking this e-process in~\eqref{eq:price} gives $\Gamma_T\le\max_{0\le t\le T}R_t(E_t)=\log\eta_T$.
\end{proof}

\subsection{Proofs of Lemma~\ref{lem:duality} and Remark~\ref{rem:deterministic}}\label{sec:proof-duality}

The proof of Lemma~\ref{lem:duality} follows the Donsker--Varadhan variational formula.
We apply Jensen's inequality to the e-process divided by the density of $Q^\xi$ with respect to $\bar P$.
The second inequality in~\eqref{eq:eproc}, applied to each stopped null distribution in $\bar P$, then gives an expectation at most one for the e-process under $\bar P$.

\begin{proof}[of Lemma~\ref{lem:duality}]
  We take $\cX^0$ to be a singleton and equip the set of the pairs $(t,x^t)$ of Section~\ref{sec:duality} with the $\sigma$-algebra generated by the sets $\{t\}\times B$, with $B\subseteq\cX^t$ measurable.
  We may assume that $\KL(Q^\xi\|\bar P)<\infty$, and hence that $Q^\xi$ has a density $f=\dd Q^\xi/\dd\bar P$.
  Since
  \[
    \int\log^-f\dd Q^\xi\le\int f^{-1}\dd Q^\xi\le1,
  \]
  the function $\log f$ is $Q^\xi$-integrable.
  Let $H(t,x^t)=E_t(x^t)$ on the set of pairs $(t,x^t)$ of Section~\ref{sec:duality}, including $t=\infty$.

  Jensen's inequality applied to the concave logarithm gives the inequality of the Donsker--Varadhan variational formula~\citep[equation~(2.1) and Theorem~2.1]{donsker1983asymptotic}
  \begin{equation}\label{eq:S-gibbs}
    \int\log H\dd Q^\xi-\KL(Q^\xi\|\bar P)=\int\log(H/f)\dd Q^\xi\le\log\int(H/f)\dd Q^\xi\le\log\int H\dd\bar P.
  \end{equation}
  By Fubini's theorem and the second inequality in~\eqref{eq:eproc},
  \begin{equation}\label{eq:S-fubini}
    \int H\dd\bar P=\int\biggl[\int_0^1\E_{P_\lambda}[E_{\tau_\lambda(X,u)}]\dd u\biggr]\pi(\dd\lambda)\le1.
  \end{equation}
  Since $\int(H/f)\dd Q^\xi\le\int H\dd\bar P\le1$, $\log^+(H/f)$ is $Q^\xi$-integrable, and the integrals in~\eqref{eq:S-gibbs} are well defined.
  Since $\int\log H\dd Q^\xi=\sum_t\xi_t\E_Q[\log E_t]$, \eqref{eq:S-gibbs} and~\eqref{eq:S-fubini} give
  \[
    \sum_t\xi_t\E_Q[\log E_t]\le\KL(Q^\xi\|\bar P).
  \]
  Lemma~\ref{lem:regret} then gives~\eqref{eq:duality}.
  For the second claim, we use $\max_{t\le T}R_t(E_t)\ge\sum_t\xi_tR_t(E_t)$ and take the infimum over $(E_t)$ and the supremum over $\bar P$.
\end{proof}

\begin{proof}[of Remark~\ref{rem:deterministic}]
  Suppose that every $\tau_\lambda(x,u)$ depends on $u$ only, and let $H(t,x^t)=\Es_t(x^t)$ for $t\le T$ and $H(t,x^t)=0$ otherwise.
  For every $\lambda$ and $u$, the stopping time $\tau_\lambda(\cdot,u)$ is then a constant $t$.
  The integral of $H$ with respect to $P_\lambda^{\tau_\lambda(\cdot,u)}$ is thus $\E_{P_\lambda}[\Es_t]\le1$ if $t\le T$, and zero otherwise.
  Integrating over $u$ and $\lambda$ gives $\int H\dd\bar P\le1$.
  Since $\int\log H\dd Q^\xi=\sum_t\xi_t\E_Q[\log\Es_t]$, \eqref{eq:S-gibbs} gives $\KL(Q^\xi\|\bar P)\ge\sum_t\xi_t\E_Q[\log\Es_t]$.
\end{proof}

\subsection{Proof of Proposition~\ref{prop:numeraire}}\label{sec:proof-numeraire}

\begin{proof}[of Proposition~\ref{prop:numeraire}]
  Let $\Ps_t$ be the mixture in the statement.
  Since $\rho^2\ge\max(0,1-s^2)$ and $t\ge1$, its mixing variance $\rho^2-(1-s^2)/t$ is nonnegative.
  We integrate~\eqref{eq:S-factor} with $v=s^2$ over $m\sim\N(0,\rho^2)$, and with $v=1$ over the mixing distribution of $\Ps_t$.
  The densities of $Q_t$ and $\Ps_t$ are then
  \[
    t^{-1/2}\varphi_{\rho^2+s^2/t}(\bar x_t)\prod_{i=2}^t\varphi_{s^2}(y_i),\quad t^{-1/2}\varphi_{\rho^2+s^2/t}(\bar x_t)\prod_{i=2}^t\varphi_1(y_i).
  \]
  Thus, $E=\dd Q_t/\dd\Ps_t=\prod_{i=2}^t\ell(Y_i)$.
  Under $P_\theta$, the $Y_i$ are independent $\N(0,1)$ variables by Lemma~\ref{lem:orth:iii}.
  For $Z\sim\N(0,1)$, $\E[\ell(Z)]=\int\varphi_{s^2}=1$.
  The products $\prod_{i=2}^t\ell(Y_i)$, $t\ge0$, therefore form a test martingale with respect to $(\cG_t)$.
  In particular, $E\in\cE_t$.
  For every $E'\in\cE_t$, Fubini's theorem gives
  \[
    \E_Q[E'/E]=\E_{\Ps_t}[E']=\int\E_\theta[E']\N\{0,\rho^2-(1-s^2)/t\}(\dd\theta)\le1.
  \]
  Consequently, $E=\Es_t$, and $\Ps_t$ is the measure in~\eqref{eq:numeraire}.
  Under every $Q_\mu$, and hence under $Q$, the $Y_i$ are independent $\N(0,s^2)$ variables.
  Therefore, $\E_Q[\log\Es_t]=(t-1)\kappa$.
\end{proof}

\subsection{Proof of Lemma~\ref{lem:identity}}\label{sec:proof-identity}

\begin{proof}[of Lemma~\ref{lem:identity}]
  (i) Equation~\eqref{eq:S-factor} with $m=\theta$ and $v\in\{s^2,1\}$ gives
  \[
    L^\theta_t=\frac{\varphi_{s^2/t}(\bar X_t-\theta)}{\varphi_{1/t}(\bar X_t-\theta)}\prod_{i=2}^t\ell(Y_i)=\frac{1}{s}\exp\biggl\{-\frac{(1-s^2)t(\bar X_t-\theta)^2}{2s^2}\biggr\}\Es_t.
  \]

  (ii) Fix $\theta\in\R$ and $T\ge0$, and let $g_0=1$ and $g_t=s\exp\{(1-s^2)t(\bar X_t-\theta)^2/(2s^2)\}$ for $t\ge1$.
  By part~\ref{lem:identity:i*}, $\Es_t=L^\theta_tg_t$ for every $t\ge0$.
  Let $\tau$ be a stopping time of $(\cF_t)$.
  Since $L^\theta_t$ is the density of $Q_\theta$ with respect to $P_\theta$ on $\cF_t$, and $g_t\ind(\tau\wedge T=t)$ is $\cF_t$-measurable,
  \[
    \E_\theta[\Es_{\tau\wedge T}]=\sum_{t=0}^T\E_\theta[L^\theta_tg_t\ind(\tau\wedge T=t)]=\sum_{t=0}^T\E_{Q_\theta}[g_t\ind(\tau\wedge T=t)]=\E_{Q_\theta}[g_{\tau\wedge T}].
  \]
  Let $S_0=0$ and $S_t=\sum_{i=1}^t(X_i-\theta)/s$ for $t\ge1$.
  Under $Q_\theta$, $(S_t)$ is a random walk with independent $\N(0,1)$ increments, as in~\ref{lem:identity:ii*}.
  Since $t(\bar X_t-\theta)^2=s^2S_t^2/t$, we have $g_t=h_t$ for every $t\ge0$.
  The map that sends the data to $(S_t)_{t\ge1}$ sends the data filtration onto the filtration generated by $(S_t)$.
  The stopping times $\tau$ of $(\cF_t)$ therefore correspond bijectively to the stopping times $\hat\tau$ of this filtration, with
  \[
    \E_\theta[\Es_{\tau\wedge T}]=\E[h_{\hat\tau\wedge T}].
  \]
  Taking the supremum over $\tau$, and then over $\theta$ in~\eqref{eq:snell}, gives the claim.
\end{proof}

\subsection{Proof of Proposition~\ref{prop:mixture}}\label{sec:proof-mixture}

The upper bound in Proposition~\ref{prop:mixture:i} follows from Proposition~\ref{prop:basic:iii} and an upper bound on $\eta_T$.
By Lemma~\ref{lem:identity:ii}, $\eta_T$ is at most the expectation of $\max_{t\le T}h_t$.
For $q>1$ with $(1-s^2)q<1$, we split the times $t\le T$ into the $O(\log T)$ blocks with $q^k\le t<q^{k+1}$.
In each block, the probability that $S_t^2/t\ge z^2$ for some $t$ is at most $2\exp\{-z^2/(2q)\}$ by Doob's maximal inequality.
A union bound over the blocks gives Lemma~\ref{lem:upper}, and hence $\log\eta_T\le(1-s^2)q\log\log T+O(1)$.

For the lower bound, we apply Lemma~\ref{lem:duality} to a mixture $\bar P$ of null distributions $P_\theta$ with $\theta\sim\N(0,\rho^2+s^2)$.
Each $P_\theta$ is stopped at the first horizon $t_k=r^k$ at which $t_k^{1/2}|\bar X_{t_k}-\theta|$ lies in an interval $[c,c+1]$, or, with probability $1/2$, at an independent uniform horizon $t_k$.
In Lemma~\ref{lem:lbmix}, we compute the density of $Q^\xi$ with respect to $\bar P$ and reduce Lemma~\ref{lem:duality} to a lower bound in terms of functions $J_k$.
Lemma~\ref{lem:jk} gives lower bounds on $J_k$.
On an event $\{Y^{t_k}\in A_k\}$, the rule stops at $t_k$ whenever $t_k^{1/2}|\bar X_{t_k}-\theta|$ lies in $[c,c+1]$, and $\log\{(K+1)J_k\}$ is then at least about $\log K-c^2/2$.
In this bound, the uniform weights on the $K+1$ horizons give the term $\log K$, and the probability of the interval $[c,c+1]$ under $\N(0,1)$ gives the term $-c^2/2$.
On the event $\{Y^{t_k}\in A_k\}$, the increments $t_j^{1/2}(\bar X_{t_k}-\bar X_{t_j})$ of the running mean stay below about $c$ in absolute value for every $j<k$.
Under $Q$, each of these increments has variance at most $s^2$.
In Lemma~\ref{lem:fixedk}, a union bound over the $K$ earlier horizons therefore shows that the event $\{Y^{t_k}\in A_k\}$ has probability close to one under $Q$ when $c$ is at least about $s(2\log K)^{1/2}$.
With this level, $\log K-c^2/2$ is close to $(1-s^2)\log K$.
Since $K$ is of order $\log T/\log r$, the lower bound divided by $\log\log T$ tends to $1-s^2$ as $T\to\infty$ and then $r\to\infty$.

The terms of the e-process of Proposition~\ref{prop:mixture:ii} are the numéraires multiplied by the deterministic factor $\sum_{j:T_j\ge t}1/\{j(j+1)u_j\}$.
It is a mixture of the e-processes of Proposition~\ref{prop:basic:iii} at the horizons $T_j$, with $\eta_{T_j}$ replaced by the bounds $u_j\ge\eta_{T_j}$ of Lemma~\ref{lem:upper}.
At a horizon $t$ with $T_{j-1}<t\le T_j$, its regret is at most $\log u_j+\log\{j(j+1)\}$.
This bound equals $(1-s^2)\log\log t$ up to terms of order $\log\log\log t$.

The following lemma gives an explicit upper bound on $\eta_T$.

\begin{lemma}[Upper bound on $\eta_T$]\label{lem:upper}
  Let $s\in(0,1)$, $a=1-s^2$, $T\ge1$, and $q\in(1,1/a)$.
  Then
  \[
    \log\eta_T\le aq\log(2+2\log T/\log q)+\log\{1/(1-aq)\}.
  \]
\end{lemma}

\begin{proof}[of Lemma~\ref{lem:upper}]
  Let $M=\max[1,\max_{1\le t\le T}\exp\{aS_t^2/(2t)\}]$, with $S_t$ as in Lemma~\ref{lem:identity:ii}.
  Since $s<1$, $h_t\le M$ for $0\le t\le T$.
  Lemma~\ref{lem:identity:ii} then gives $\eta_T\le\E[M]$.

  For $k\ge0$, let $B_k=\{t\le T:q^k\le t<q^{k+1}\}$.
  Since $B_k$ is empty for $k>\log T/\log q$, at most $N=1+\lfloor\log T/\log q\rfloor$ blocks are nonempty.
  If $B_k$ is nonempty, its smallest element $m$ and its largest element $u$ satisfy $u/m<q$.
  For $x>0$, Doob's maximal inequality applied to the submartingale $\exp(xS_t/u)$ gives
  \begin{equation}\label{eq:S-doob}
    \pr\Bigl(\max_{t\le u}S_t\ge x\Bigr)\le\exp\{-x^2/(2u)\},
  \end{equation}
  where $\pr$ denotes probability.
  By symmetry,~\eqref{eq:S-doob} holds with $|S_t|$ in place of $S_t$ and an additional factor $2$.
  Therefore, for $z>0$,
  \begin{equation}\label{eq:S-block}
    \pr\Bigl(\max_{t\in B_k}S_t^2/t\ge z^2\Bigr)\le\pr\Bigl(\max_{t\le u}|S_t|\ge zm^{1/2}\Bigr)\le2\exp\{-z^2m/(2u)\}\le2\exp\{-z^2/(2q)\}.
  \end{equation}
  A union bound in~\eqref{eq:S-block} over the nonempty blocks gives
  \[
    \pr\Bigl(\max_{1\le t\le T}S_t^2/t\ge z^2\Bigr)\le2N\exp\{-z^2/(2q)\}.
  \]

  For $x\ge1$, $M>x$ exactly when $\max_{1\le t\le T}S_t^2/t>(2/a)\log x$.
  Hence,
  \[
    \pr(M>x)\le\min\{1,2Nx^{-1/(aq)}\}.
  \]
  With $p=1/(aq)>1$ and $N'=2N\ge1$,
  \[
    \E[M]=1+\int_1^\infty\pr(M>x)\dd x\le1+\int_1^\infty\min(1,N'x^{-p})\dd x=\frac{pN'^{1/p}}{p-1}=\frac{N'^{aq}}{1-aq}.
  \]
  Since $N'\le2+2\log T/\log q$, the claim follows.
\end{proof}

\begin{proof}[of Proposition~\ref{prop:mixture:i}, upper bound]
  By Proposition~\ref{prop:basic:iii}, $\Gamma_T\le\log\eta_T$.
  For fixed $q\in(1,1/a)$, we have $\log(2+2\log T/\log q)=\log\log T+O(1)$, and Lemma~\ref{lem:upper} gives
  \[
    \limsup_{T\to\infty}\frac{\log\eta_T}{\log\log T}\le aq.
  \]
  Letting $q\to1$ replaces $aq$ by $a=1-s^2$.
  Together with the lower bound on $\Gamma_T$ proved below, the inequality $\Gamma_T\le\log\eta_T$ gives both limits in~\ref{prop:mixture:i*}.
\end{proof}

For the lower bound, let $s\in(0,1)$ and $\rho^2\ge1-s^2$, as in Proposition~\ref{prop:mixture}, and fix integers $r\ge4$ and $K\ge1$, a horizon $T\ge r^K$, and a constant $c>0$.
Let $t_k=r^k$ for $0\le k\le K$, let $\xi\in\Delta_T$ be uniform on $\{t_0,\ldots,t_K\}$, and let $\gamma^2=\rho^2+s^2$.
For $\theta\in\R$, let $Z_t=t^{1/2}(\bar X_t-\theta)$, and let $\tau'$ be the first $t_k$ with $c\le|Z_{t_k}|\le c+1$, or $\infty$ if there is none.
Let
\[
  \bar P=\int\biggl(\frac12P_\theta^{\tau'}+\frac12\sum_{k=0}^K\xi_{t_k}P_\theta^{t_k}\biggr)\N(0,\gamma^2)(\dd\theta),
\]
where $P_\theta^{t_k}$ is $P_\theta$ stopped at the constant time $t_k$.

For $0\le k\le K$, write $t=t_k$, and write $\bar x_t$ and $y=(y_2,\ldots,y_t)$ for the functions of $x\in\R^t$ in Lemma~\ref{lem:orth}.
For $j<k$, let
\[
  W_{jk}(y)=t_j^{1/2}\sum_{i=t_j+1}^{t}\frac{y_i}{\{i(i-1)\}^{1/2}},
\]
and let
\[
  \mathcal Z_k(y)=\{z\in\R:c\le|z|\le c+1\}\setminus\bigcup_{j<k}\{z\in\R:c\le|(t_j/t)^{1/2}z-W_{jk}(y)|\le c+1\}.
\]
Let
\begin{equation}\label{eq:S-J}
  J_k(x)=\frac{1}{2\varphi_{\rho^2+s^2/t}(\bar x_t)}\biggl\{\int_{\mathcal Z_k(y)}\varphi_1(z)\varphi_{\gamma^2}(\bar x_t-t^{-1/2}z)\dd z+\xi_t\varphi_{\gamma^2+1/t}(\bar x_t)\biggr\},
\end{equation}
and write $J_k$ for $J_k(X^{t_k})$.

\begin{lemma}[Mixture of stopped null distributions]\label{lem:lbmix}
  The distribution $\bar P$ belongs to $\cB$, and
  \begin{equation}\label{eq:S-lb}
    \Gamma_T\ge\frac{1}{K+1}\sum_{k=0}^K\E_Q[\log\{(K+1)J_k\}].
  \end{equation}
\end{lemma}

\begin{proof}[of Lemma~\ref{lem:lbmix}]
  For each $\theta$, let $\tau_\theta(\cdot,u)=\tau'$ for $u<1/2$, and let $\tau_\theta(\cdot,u)=t_k$ for $u\ge1/2$, where $k$ is the integer part of $(K+1)(2u-1)$, or $K$ if $u=1$.
  Since $\xi_{t_k}=1/(K+1)$, the randomized stopping time $\tau_\theta$ gives
  \[
    P_\theta^{\tau_\theta}=\int_0^1P_\theta^{\tau_\theta(\cdot,u)}\dd u=\frac12P_\theta^{\tau'}+\frac12\sum_{k=0}^K\xi_{t_k}P_\theta^{t_k}.
  \]
  Hence, $\bar P=\int P_\theta^{\tau_\theta}\N(0,\gamma^2)(\dd\theta)$ belongs to $\cB$.

  Fix $k$, write $t=t_k$, and write densities with respect to the Lebesgue measure on $\R^t$ as functions of $\bar x_t$ and $y$.
  By the proof of Proposition~\ref{prop:numeraire}, the restriction of $Q^\xi$ to $\{t\}\times\R^t$ has density
  \[
    \xi_tt^{-1/2}\varphi_{\rho^2+s^2/t}(\bar x_t)\prod_{i=2}^t\varphi_{s^2}(y_i).
  \]
  Under $P_\theta$,~\eqref{eq:S-factor} gives $X^t$ the density
  \[
    t^{-1/2}\varphi_{1/t}(\bar x_t-\theta)\prod_{i=2}^t\varphi_1(y_i).
  \]
  By Lemma~\ref{lem:orth:ii}, $Z_{t_j}=(t_j/t)^{1/2}Z_t-W_{jk}(Y^t)$ for $j<k$.
  Thus, $\tau'=t$ exactly when $Z_t\in\mathcal Z_k(Y^t)$.
  After the substitution $\theta=\bar x_t-t^{-1/2}z$ in the term with $P_\theta^{\tau'}$, the restriction of $\bar P$ to $\{t\}\times\R^t$ has density
  \[
    t^{-1/2}\varphi_{\rho^2+s^2/t}(\bar x_t)\prod_{i=2}^t\varphi_1(y_i)J_k(x).
  \]
  Hence, $\dd Q^\xi/\dd\bar P=\xi_t\Es_t/J_k$ on $\{t\}\times\R^t$.

  Taking the expectation of the logarithm of this density under $Q^\xi$ gives
  \[
    \KL(Q^\xi\|\bar P)=\sum_{k=0}^K\xi_{t_k}\E_Q[\log\Es_{t_k}]-\sum_{k=0}^K\xi_{t_k}\E_Q[\log(J_k/\xi_{t_k})].
  \]
  Lemma~\ref{lem:duality} with $\xi_{t_k}=1/(K+1)$ therefore gives~\eqref{eq:S-lb}.
\end{proof}

\begin{lemma}[Lower bounds on $J_k$]\label{lem:jk}
  Let $0\le k\le K$ and $t=t_k$, let $\beta'=\tfrac12\log\{\rho^2/(\gamma^2+1)\}-\log2$, and let
  \[
    p_k(y)=\int_{\mathcal Z_k(y)}\varphi_1(z)\dd z.
  \]
  \begin{parts}
  \item\partlabel{lem:jk:i} For every $x\in\R^t$, $\log\{(K+1)J_k(x)\}\ge\beta'$.
  \item\partlabel{lem:jk:ii} For every $x\in\R^t$,
    \[
      J_k(x)\ge\frac{\rho}{2\gamma}\exp\biggl\{-\frac{2(c+1)t^{-1/2}|\bar x_t|+(c+1)^2t^{-1}}{2\gamma^2}\biggr\}p_k(y).
    \]
  \item\partlabel{lem:jk:iii} Let $c'=c-(c+1)r^{-1/2}$, assume that $c'>0$, and let $A_k=\{y:\max_{j<k}|W_{jk}(y)|<c'\}$.
    Then $p_k(y)\ge2\varphi_1(c+1)$ for every $y\in A_k$.
  \end{parts}
\end{lemma}

Since $Z_t\sim\N(0,1)$ is independent of $Y^t$ under $P_\theta$, $p_k(Y^t)$ is the conditional probability that $\tau'=t$, given $Y^t$.
The proof of part~\ref{lem:jk:iii*} shows that, on the event $\{Y^t\in A_k\}$, the rule $\tau'$ stops at $t$ whenever $c\le|Z_t|\le c+1$.

\begin{proof}[of Lemma~\ref{lem:jk}]
  (i) For $0<v_1\le v_2$, $\varphi_{v_2}/\varphi_{v_1}\ge(v_1/v_2)^{1/2}$ pointwise.
  Since $\gamma^2+1/t\ge\rho^2+s^2/t\ge\rho^2$ and $\xi_t=1/(K+1)$, the second term of $J_k$ in~\eqref{eq:S-J} gives
  \[
    (K+1)J_k(x)\ge\frac{\varphi_{\gamma^2+1/t}(\bar x_t)}{2\varphi_{\rho^2+s^2/t}(\bar x_t)}\ge\frac12\biggl(\frac{\rho^2}{\gamma^2+1}\biggr)^{1/2}.
  \]
  Taking logarithms gives part~\ref{lem:jk:i*}.

  (ii) Since $|z|\le c+1$ for $z\in\mathcal Z_k(y)$,
  \[
    (\bar x_t-t^{-1/2}z)^2\le\bar x_t^2+2(c+1)t^{-1/2}|\bar x_t|+(c+1)^2t^{-1}.
  \]
  Hence, for $z\in\mathcal Z_k(y)$,
  \[
    \varphi_{\gamma^2}(\bar x_t-t^{-1/2}z)\ge\varphi_{\gamma^2}(\bar x_t)\exp\biggl\{-\frac{2(c+1)t^{-1/2}|\bar x_t|+(c+1)^2t^{-1}}{2\gamma^2}\biggr\}.
  \]
  Since $\rho^2+s^2/t\le\gamma^2$, the inequality in the proof of part~\ref{lem:jk:i*} gives $\varphi_{\gamma^2}(\bar x_t)/\varphi_{\rho^2+s^2/t}(\bar x_t)\ge\rho/\gamma$.
  The first term of $J_k$ in~\eqref{eq:S-J} therefore gives part~\ref{lem:jk:ii*}.

  (iii) Let $y\in A_k$, $j<k$, and $c\le|z|\le c+1$.
  Since $(t_j/t)^{1/2}=r^{-(k-j)/2}\le r^{-1/2}$,
  \[
    |(t_j/t)^{1/2}z-W_{jk}(y)|<(c+1)r^{-1/2}+c'=c.
  \]
  Hence, $\mathcal Z_k(y)=\{z:c\le|z|\le c+1\}$, and
  \[
    p_k(y)=2\int_c^{c+1}\varphi_1(z)\dd z\ge2\varphi_1(c+1).
  \]
\end{proof}

\begin{lemma}[Lower bound for fixed $K$]\label{lem:fixedk}
  Let $c'$, $A_k$, and $\beta'$ be as in Lemma~\ref{lem:jk}, and let $\beta=\log(K+1)+\log\varphi_1(c+1)+\log(\rho/\gamma)$.
  If $c'>0$ and $\beta\ge\beta'$, then
  \begin{equation}\label{eq:S-final}
    \Gamma_T\ge\beta-(\beta-\beta')K\exp\{-c'^2/(2s^2)\}-\frac{2(c+1)/\gamma+2(c+1)^2/(3\gamma^2)}{K+1}.
  \end{equation}
\end{lemma}

\begin{proof}[of Lemma~\ref{lem:fixedk}]
  Fix $k$ and write $t=t_k$.
  Under $Q$, we have $W_{jk}(Y^t)\sim\N\{0,s^2(1-t_j/t)\}$ for $j<k$.
  A $\N(0,v)$ variable has absolute value at least $x\ge0$ with probability at most $\exp\{-x^2/(2v)\}$.
  A union bound over $j<k$ therefore gives
  \[
    Q(Y^t\notin A_k)\le K\exp\{-c'^2/(2s^2)\}.
  \]

  On the event $\{Y^t\in A_k\}$, Lemma~\ref{lem:jk:ii} and~\ref{lem:jk:iii*} give
  \[
    \log\{(K+1)J_k\}\ge\beta-\frac{2(c+1)t^{-1/2}|\bar X_t|+(c+1)^2t^{-1}}{2\gamma^2}.
  \]
  On its complement, Lemma~\ref{lem:jk:i} gives $\log\{(K+1)J_k\}\ge\beta'$.
  Since $\beta'=\beta-(\beta-\beta')$, both bounds give, everywhere,
  \[
    \log\{(K+1)J_k\}\ge\beta-(\beta-\beta')\ind(Y^t\notin A_k)-\frac{2(c+1)t^{-1/2}|\bar X_t|+(c+1)^2t^{-1}}{2\gamma^2}.
  \]
  Under $Q$, $\bar X_t\sim\N(0,\rho^2+s^2/t)$, and hence $\E_Q[|\bar X_t|]\le(\rho^2+s^2/t)^{1/2}\le\gamma$.
  Since $\beta-\beta'\ge0$, taking expectations under $Q$ gives
  \[
    \E_Q[\log\{(K+1)J_k\}]\ge\beta-(\beta-\beta')K\exp\{-c'^2/(2s^2)\}-\frac{c+1}{\gamma t_k^{1/2}}-\frac{(c+1)^2}{2\gamma^2t_k}.
  \]

  Since $r\ge4$, we have $\sum_kr^{-k/2}\le2$ and $\sum_kr^{-k}\le4/3$.
  Averaging over $k$ in~\eqref{eq:S-lb} therefore gives~\eqref{eq:S-final}.
\end{proof}

\begin{proof}[of Proposition~\ref{prop:mixture:i}, lower bound]
  For $K\ge3$, choose
  \[
    c'=s(2\log K+4\log\log K)^{1/2},\quad c=(c'+r^{-1/2})/(1-r^{-1/2}),
  \]
  so that $c-(c+1)r^{-1/2}=c'>0$ and $K\exp\{-c'^2/(2s^2)\}=(\log K)^{-2}$.
  Then $c+1=(c'+1)/(1-r^{-1/2})=O\{(\log K)^{1/2}\}$, and
  \[
    (c+1)^2=\frac{2s^2\log K}{(1-r^{-1/2})^2}+O\{(\log K)^{1/2}\}.
  \]
  Since $\log\varphi_1(c+1)=-(c+1)^2/2+O(1)$,
  \[
    \beta=\{1-s^2(1-r^{-1/2})^{-2}\}\log K+O\{(\log K)^{1/2}\}.
  \]
  By the choice of $c'$ and since $c+1=O\{(\log K)^{1/2}\}$, the last two terms of~\eqref{eq:S-final} tend to zero as $K\to\infty$.
  Since $\Gamma_T\ge0$, the bound below holds trivially if the coefficient of $\log K$ in $\beta$ is not positive.
  Otherwise, $\beta\ge\beta'$ for large $K$, and Lemma~\ref{lem:fixedk} applies.
  For $T\ge r^3$, the choice $K=\lfloor\log T/\log r\rfloor$ satisfies $r^K\le T$ and $\log K=\log\log T-\log\log r+o(1)$.
  Then~\eqref{eq:S-final} gives
  \[
    \liminf_{T\to\infty}\frac{\Gamma_T}{\log\log T}\ge1-\frac{s^2}{(1-r^{-1/2})^2},
  \]
  and letting $r\to\infty$ gives the lower bound in~\ref{prop:mixture:i*}.
\end{proof}

\begin{proof}[of Proposition~\ref{prop:mixture:ii}]
  Since $(1-s^2)q_j=1-s^2j/(j+1)$, we have $1<q_j<1/(1-s^2)$, and Lemma~\ref{lem:upper} with $T=T_j$ and $q=q_j$ gives $\eta_{T_j}\le u_j$.
  The weights $c_j=1/\{j(j+1)u_j\}$ therefore satisfy $\sum_jc_j\eta_{T_j}\le\sum_j1/\{j(j+1)\}=1$.
  Let $\theta\in\R$ and let $\tau$ be a stopping time of $(\cF_t)$.
  Since $\Es_t\ind(t\le T_j)\le\Es_{t\wedge T_j}$ for all $t$ and $j$, the definition of $\eta_{T_j}$ in~\eqref{eq:snell} and monotone convergence give, for every $n\ge0$,
  \[
    \E_\theta[E_{\tau\wedge n}]=\sum_jc_j\E_\theta[\Es_{\tau\wedge n}\ind(\tau\wedge n\le T_j)]\le\sum_jc_j\E_\theta[\Es_{\tau\wedge n\wedge T_j}]\le\sum_jc_j\eta_{T_j}\le1.
  \]
  Under $P_\theta$, the $Y_i$ are independent $\N(0,1)$ variables with $\E_\theta[\log\ell(Y_i)]=-\KL\{\N(0,1)\|\N(0,s^2)\}<0$.
  Hence, $\Es_t\to0$ almost surely.
  Since $E_t\le\Es_t\sum_jc_j\le\Es_t$, also $E_t\to0$ almost surely.
  Since $E_{\tau\wedge n}\to E_\tau$ almost surely as $n\to\infty$, Fatou's lemma gives $\E_\theta[E_\tau]\le1$.

  For the upper limit, Lemma~\ref{lem:regret} and the definition of $E_t$ give $R_t(E_t)=-\log\sum_{j:T_j\ge t}c_j$.
  For $j\ge1$ and $T_{j-1}<t\le T_j$, with $T_0=0$, this sum is at least $c_j$, and hence $R_t(E_t)\le\log u_j+\log\{j(j+1)\}$.
  Since $\log T_j=e^j+O(1)$ and $1/\log q_j=O(j)$, we have $\log u_j=(1-s^2)j+O(\log j)$.
  Since $\log\log t>\log\log T_{j-1}=j-1+o(1)$, we have $j\le\log\log t+1+o(1)$ as $t\to\infty$.
  Consequently, $\limsup_{t\to\infty}R_t(E_t)/\log\log t\le1-s^2$.

  For the lower limit, let $t_T$ maximize $R_t(E_t)$ over $t\le T$.
  The previous paragraph gives $R_t(E_t)<\infty$ for every $t$.
  Since $R_{t_T}(E_{t_T})\ge\Gamma_T\to\infty$ by part~\ref{prop:mixture:i*}, it follows that $t_T\to\infty$.
  Since $t_T\le T$,
  \[
    \frac{R_{t_T}(E_{t_T})}{\log\log t_T}\ge\frac{\Gamma_T}{\log\log T},
  \]
  and part~\ref{prop:mixture:i*} gives $\limsup_{t\to\infty}R_t(E_t)/\log\log t\ge1-s^2$.
\end{proof}

\begin{remark}\label{rem:stopped}
  In the first version of this paper, the e-process of Proposition~\ref{prop:mixture:ii} was the mixture $(E'_t)_{t\ge0}$ of the numéraires stopped at the horizons $T_j$, with
  \[
    E'_t=\sum_{j=1}^\infty\frac{\Es_{t\wedge T_j}}{j(j+1)u_j}.
  \]
  The bounds $\eta_{T_j}\le u_j$ and Fatou's lemma also make $(E'_t)$ an e-process.
  Since $\Es_t\ind(t\le T_j)\le\Es_{t\wedge T_j}$, the process $(E'_t)$ dominates the e-process $(E_t)$ of Proposition~\ref{prop:mixture:ii}.
  Its regrets are therefore at most those of $(E_t)$, and the last paragraph of the proof of Proposition~\ref{prop:mixture:ii} shows that they have the same limit.
  However, for $t>T_j$, the stopped numéraire $\Es_{T_j}$ in $E'_t$ depends only on $X^{T_j}$.
  Following a question of Peter Gr\"unwald on e-processes with terms invariant under permutations of the observations, we replaced each $\Es_{t\wedge T_j}$ by $\Es_t\ind(t\le T_j)$.
  The terms of the resulting e-process $(E_t)$ are invariant under permutations of $X_1,\ldots,X_t$ (Section~\ref{sec:rate}), and its regret is $-\log\sum_{j:T_j\ge t}c_j$.
\end{remark}

\subsection{Proof of Theorem~\ref{thm:composite}}\label{sec:proof-composite}

For the lower bound in Theorem~\ref{thm:composite:i}, we average the regret against $Q_\mu$ over $\mu\sim\N(0,\rho^2)$.
Since the numéraires of $Q_\mu$ and of $Q$ have the e-powers $t\kappa$ and $(t-1)\kappa$, this average equals $\kappa$ plus the regret against $Q$.
For the upper bound and part~\ref{thm:composite:ii*}, the e-processes of Propositions~\ref{prop:basic:iii} and~\ref{prop:mixture:ii} are functions of the $Y_i$, and the distribution of the $Y_i$ is the same under every $Q_\mu$ and under $Q$.

\begin{proof}[of Theorem~\ref{thm:composite}]
  By Proposition~\ref{prop:simple:i}, the numéraire of $Q_\mu$ has e-power $t\kappa$.
  Lemma~\ref{lem:regret} then gives $R^\mu_t(E_t)=t\kappa-\E_{Q_\mu}[\log E_t]$.

  (i) We first prove the lower bound, for a fixed $T\ge1$.
  Since $\cF_0$ is trivial, every stopping time $\tau$ of $(\cF_t)$ satisfies $\tau=0$ everywhere or $\tau\ge1$ everywhere.
  Replacing $E_0$ by $1$ therefore preserves the e-process property and increases no regret.
  We may thus restrict both infima to e-processes with $E_0=1$.
  For these e-processes, the regret at $t=0$ vanishes under every alternative.
  Let $(E_t)$ be such an e-process and $1\le t\le T$.
  Since $Q=\int Q_\mu\N(0,\rho^2)(\dd\mu)$ and $\E_Q[\log\Es_t]=(t-1)\kappa$ by Proposition~\ref{prop:numeraire},
  \[
    \sup_{\mu\in\R}R^\mu_t(E_t)\ge\int\{t\kappa-\E_{Q_\mu}[\log E_t]\}\N(0,\rho^2)(\dd\mu)=t\kappa-\E_Q[\log E_t]=\kappa+R_t(E_t),
  \]
  where the exchange of integrals holds separately for the positive and negative parts of $\log E_t$.
  Since the regret at $t=0$ vanishes, the maximum over $1\le t\le T$ equals the maximum over $0\le t\le T$, both in $\Gamma_T$ and in the middle term of Theorem~\ref{thm:composite}.
  Taking this maximum and then the infimum over $(E_t)$ gives the lower bound.

  For the upper bound, we use the e-process $(E_t)_{t\ge0}$ with $E_t=\Es_t\ind(t\le T)/\eta_T$ of Proposition~\ref{prop:basic:iii}.
  Since the distributions of $(Y_i)_{i\ge2}$ under $Q_\mu$ and under $Q$ coincide, $\E_{Q_\mu}[\log\Es_t]=(t-1)\kappa$.
  For $1\le t\le T$ and every $\mu$, $R^\mu_t(E_t)=\kappa+\log\eta_T$.
  At $t=0$, $R^\mu_0(E_0)=\log\eta_T$.
  The computation of the lower bound also gives, for $t\ge1$ and every $E\in\cE_t$,
  \[
    \sup_\mu R^\mu_t(E)\ge\kappa+R_t(E)\ge\kappa,
  \]
  with equality for $E=\Es_t$.
  Proposition~\ref{prop:mixture:i} gives the limit.

  (ii) Each term $E_t$ of the e-process of Proposition~\ref{prop:mixture:ii} is a function of $Y^t$.
  The distribution of $Y^t$ is the same under every $Q_\mu$ and under $Q$.
  Hence, $\E_{Q_\mu}[\log E_t]=\E_Q[\log E_t]$ for every $\mu$.
  Since the numéraires of $Q_\mu$ and of $Q$ have the e-powers $t\kappa$ and $(t-1)\kappa$, Lemma~\ref{lem:regret} gives, for $t\ge1$,
  \[
    R^\mu_t(E_t)=t\kappa-\E_Q[\log E_t]=\kappa+R_t(E_t).
  \]
  Proposition~\ref{prop:mixture:ii} gives the limit.
\end{proof}

\subsection{Proof of Proposition~\ref{prop:sge1}}\label{sec:proof-sge1}

\begin{proof}[of Proposition~\ref{prop:sge1}]
  If $s\ge1$, then $h_t=s\exp\{-(s^2-1)S_t^2/(2t)\}\le s$ for $t\ge1$ and $h_0=1\le s$.
  It follows from Proposition~\ref{prop:basic:iii} and Lemma~\ref{lem:identity:ii} that $\Gamma_T\le\log\eta_T\le\log s$.

  For the limit, let $\varepsilon>0$ and let $\tau_\varepsilon$ be the first $t\ge1$ with $S_t^2\le\varepsilon t$.
  If $S_{t-1}S_t\le0$, then $|S_t|\le|S_t-S_{t-1}|$.
  Since the random walk $S$ changes sign infinitely often and $(S_t-S_{t-1})^2/t\to0$ almost surely, $\tau_\varepsilon<\infty$ almost surely.
  Therefore,
  \begin{equation}\label{eq:S-sge1}
    \eta_T\ge\E[h_{\tau_\varepsilon\wedge T}]\ge s\exp\{-(s^2-1)\varepsilon/2\}\pr(\tau_\varepsilon\le T),
  \end{equation}
  and the right-hand side of~\eqref{eq:S-sge1} tends to $s\exp\{-(s^2-1)\varepsilon/2\}$ as $T\to\infty$.
  Letting $\varepsilon\to0$ gives $\liminf_{T\to\infty}\eta_T\ge s$.
\end{proof}

\section{Additional results}\label{sec:additional}

\subsection{Finite sample spaces and absolute continuity}\label{sec:finite}

Both inequalities in~\eqref{eq:gamma-lb} are equalities in the following sense.
By Sion's minimax theorem~\citep[Corollary~3.3]{sion1958general}, the first is an equality for the best $\xi\in\Delta_T$ on every sample space.
On finite sample spaces, the primal-dual equivalence of \citet[Theorem~5]{koolen2022log} makes the second an equality.
The versions of Remark~\ref{rem:zero} exist when $X_1,X_2,\ldots$ are independent and identically distributed under $Q$ and under every $P\in\cP$, provided the reverse information projection of $Q_1$ on the set $\{P_1:P\in\cP\}$ of null distributions of $X_1$ belongs to this set~\citep[Theorem~12]{koolen2022log}.
This condition holds for some exponential families~\citep[Theorem~1]{grunwald2024optimal}.
The following proposition states the converse in Remark~\ref{rem:zero} and the two equalities.
This converse also holds when the null distributions of $X^T$ are absolutely continuous with respect to $Q_T$.

\begin{proposition}\label{prop:minimax}
  Let $T\ge0$.
  \begin{parts}
  \item\partlabel{prop:minimax:i} Suppose that $\Gamma_T=0$, and that $\cX$ is finite or $P_T$ is absolutely continuous with respect to $Q_T$ for every $P\in\cP$.
    Then $\Es_0,\ldots,\Es_T$ have versions such that $(\Es_{t\wedge T})_{t\ge0}$ is an e-process.
  \item\partlabel{prop:minimax:ii} On every sample space $\cX$, $\Gamma_T=\max_{\xi\in\Delta_T}\inf_{(E_t)}\sum_{t=0}^T\xi_tR_t(E_t)$.
  \item\partlabel{prop:minimax:iii} If $\cX$ is finite, then, for every $\xi\in\Delta_T$,
    \begin{equation}\label{eq:avg_dual}
      \inf_{(E_t)}\sum_{t=0}^T\xi_tR_t(E_t)=\sum_{t=0}^T\xi_t\E_Q[\log\Es_t]-\inf_{\bar P\in\cB}\KL(Q^\xi\|\bar P).
    \end{equation}
  \end{parts}
\end{proposition}

Consider the space of the pairs $(t,x^t)$ of Section~\ref{sec:duality}.
On this space, the stopped null distributions $P^\tau$ form the null, and the distributions $Q^\xi$ with $\xi$ a unit mass form the alternatives.
Then $\Gamma_T$ is the minimax pseudo-regret of \citet[Section~4]{grunwald2024safe}.
Against the alternative $Q^\xi$ with $\xi$ the unit mass at $t$, the pseudo-regret of an e-process $(E_t)$ is $R_t(E_t)$.
On a finite $\cX$, the e-processes stopped at $T$ may be taken bounded.
Indeed, if a node $x^t$ with $t\le T$ had $Q(X^t=x^t)>0$ and $P(X^t=x^t)=0$ for every $P\in\cP$, then the e-power $Q(X^t=x^t)\log(1+m)$ of the e-variable $1+m\ind(X^t=x^t)$ would tend to infinity with $m$, against the assumption of Section~\ref{sec:setting}.
Hence, every node with $Q(X^t=x^t)>0$ has $P(X^t=x^t)>0$ for some $P\in\cP$, and every e-process satisfies $E_t(x^t)\le1/P(X^t=x^t)$ there.
Setting $E_t(x^t)$ to zero at the nodes with $Q(X^t=x^t)=0$ preserves the e-process property and the regrets.
Moreover, the alternatives $Q^\xi$ with $\xi$ a unit mass have disjoint supports.
For a finite $\cX$, Proposition~\ref{prop:minimax:ii} is thus an instance of the strong duality of \citet[Theorem~3.4]{ram2026strong} for bounded e-variables.

\begin{proof}[of Proposition~\ref{prop:minimax}]
  (i) For every e-process $(E_t)$, the process $(E_{t\wedge T})_{t\ge0}$ is an e-process with the same regrets at $t\le T$.
  Since $\Gamma_T=0$, there are thus e-processes $(E^n_t)_{t\ge0}$, $n\ge1$, with $R_t(E^n_t)\le2^{-n}$ for $t\le T$.
  If $\cX$ is finite, we set $E^n_t$ and $\Es_t$ to zero at the nodes $x^t$ with $Q(X^t=x^t)=0$.
  This preserves the e-process property, the regrets, and the numéraires up to $Q$-null sets.
  For $t\le T$, let $U^n_t=E^n_t/\Es_t$.
  The first property in~\eqref{eq:numeraire} gives $\E_Q[U^n_t]\le1$, and Lemma~\ref{lem:regret} gives $\E_Q[\log U^n_t]=-R_t(E^n_t)$.
  Since $u-1-\log u\ge(u^{1/2}-1)^2$ for $u\ge0$,
  \[
    \E_Q[\{(U^n_t)^{1/2}-1\}^2]\le\E_Q[U^n_t-1-\log U^n_t]\le R_t(E^n_t)\le2^{-n}.
  \]
  Hence, for every $t\le T$, $\sum_n\{(U^n_t)^{1/2}-1\}^2<\infty$, and thus $U^n_t\to1$, on an event of $\cF_t$ with $Q$-probability one.
  If $P_T$ is absolutely continuous with respect to $Q_T$ for every $P\in\cP$, this event has probability one under every $P\in\cP$.
  Since $\E_P[\Es_t]\le1$, it follows that $E^n_t=U^n_t\Es_t\to\Es_t$ $P$-almost surely.
  If $\cX$ is finite, the event on which $U^n_t\to1$ contains every node with $Q(X^t=x^t)>0$.
  Since $E^n_t=\Es_t=0$ at the other nodes, $E^n_t\to\Es_t$ everywhere.

  Let $P\in\cP$, and let $\tau$ be a stopping time of $(\cF_t)$.
  Then $E^n_{\tau\wedge T}\to\Es_{\tau\wedge T}$ $P$-almost surely, and Fatou's lemma gives
  \[
    \E_P[\Es_{\tau\wedge T}]\le\liminf_{n\to\infty}\E_P[E^n_{\tau\wedge T}]\le1.
  \]
  The process $(\Es_{t\wedge T})_{t\ge0}$ is therefore an e-process, for every choice of versions of $\Es_0,\ldots,\Es_T$ if every $P_T$ is absolutely continuous with respect to $Q_T$, and for the versions set to zero at the nodes of $Q$-probability zero if $\cX$ is finite.

  (ii) Let $\mathcal V$ be the set of the vectors $v\in\R^{T+1}$ such that $v_t\ge R_t(E_t)$ for $t\le T$ for some e-process $(E_t)$.
  For e-processes $(E_t)$ and $(E'_t)$ and $\lambda\in[0,1]$, the process with terms $(1-\lambda)E_t+\lambda E'_t$ is an e-process.
  Since the logarithm is concave, Lemma~\ref{lem:regret} gives $R_t\{(1-\lambda)E_t+\lambda E'_t\}\le(1-\lambda)R_t(E_t)+\lambda R_t(E'_t)$ when both regrets are finite.
  The set $\mathcal V$ is therefore convex, and $\Gamma_T=\inf_{v\in\mathcal V}\max_{t\le T}v_t$.
  Sion's theorem~\citep[Corollary~3.3]{sion1958general}, applied to the bilinear function $(\xi,v)\mapsto\sum_t\xi_tv_t$ on the compact convex set $\Delta_T$ and the convex set $\mathcal V$, gives
  \begin{equation}\label{eq:S-sion}
    \Gamma_T=\inf_{v\in\mathcal V}\max_{\xi\in\Delta_T}\sum_{t=0}^T\xi_tv_t=\sup_{\xi\in\Delta_T}\inf_{v\in\mathcal V}\sum_{t=0}^T\xi_tv_t.
  \end{equation}
  For $\varepsilon\in(0,1)$ and an e-process $(E_t)$, the process with terms $(1-\varepsilon)E_t+\varepsilon$ is an e-process.
  Its terms are at least $\varepsilon$ and at least $(1-\varepsilon)E_t$.
  Its regrets at $t\le T$ are therefore finite, and they exceed those of $(E_t)$ by at most $\log\{1/(1-\varepsilon)\}$.
  Letting $\varepsilon\to0$ gives, for every $\xi\in\Delta_T$,
  \begin{equation}\label{eq:S-supinf}
    \inf_{v\in\mathcal V}\sum_{t=0}^T\xi_tv_t=\inf_{(E_t)}\sum_{t=0}^T\xi_tR_t(E_t).
  \end{equation}
  As an infimum of continuous linear functions of $\xi$, the left-hand side of~\eqref{eq:S-supinf} is upper semicontinuous on the compact set $\Delta_T$.
  The supremum in~\eqref{eq:S-sion} is therefore attained, and~\eqref{eq:S-supinf} gives~\ref{prop:minimax:ii*}.

  (iii) As shown after the statement, the infimum in~\eqref{eq:avg_dual} may be restricted to finite-valued e-processes stopped at $T$.
  Such an e-process is an anytime-valid e-value for the stopping times of $(\cF_t)$ bounded by $T$ in the sense of \citet[Definition~1]{koolen2022log}.
  Conversely, every such e-value defines an e-process stopped at $T$.
  By Lemma~\ref{lem:regret},
  \[
    \inf_{(E_t)}\sum_{t=0}^T\xi_tR_t(E_t)=\sum_{t=0}^T\xi_t\E_Q[\log\Es_t]-\sup_{(E_t)}\sum_{t=0}^T\xi_t\E_Q[\log E_t].
  \]
  \citet[Theorem~5]{koolen2022log} assume only a finite outcome space and a finite horizon, with no condition on the null.
  Applied to a data-independent stopping time with distribution $\xi$, their theorem gives
  \[
    \sup_{(E_t)}\sum_{t=0}^T\xi_t\E_Q[\log E_t]=\inf_{\bar P}\KL(Q^\xi\|\bar P),
  \]
  where $\bar P$ ranges over their Bayesian mixtures of stopped null distributions.
  Since these mixtures belong to $\cB$, the left-hand side of~\eqref{eq:avg_dual} is at most its right-hand side.
  Lemma~\ref{lem:duality} gives the reverse inequality.
\end{proof}

\subsection{A Bernoulli example}

For every $\xi\in\Delta_T$, Lemma~\ref{lem:regret} gives
\begin{equation}\label{eq:S-epower}
  \sup_{(E_t)}\sum_{t=0}^T\xi_t\E_Q[\log E_t]=\sum_{t=0}^T\xi_t\E_Q[\log\Es_t]-\inf_{(E_t)}\sum_{t=0}^T\xi_tR_t(E_t).
\end{equation}
\citet[Section~4.5]{koolen2022log} compute the left-hand side of~\eqref{eq:S-epower} for the uniform $\xi$ and $T\in\{3,6\}$ in the following example.

\begin{example}\label{ex:bernoulli}
  \citet[Section~4.5]{koolen2022log} consider the null and the alternative
  \[
    \cP=\{\mathrm{Ber}(\theta)^{\otimes\infty}:\theta=0.3,0.4,0.5,0.6,0.7\},\quad Q=\tfrac12\mathrm{Ber}(0.2)^{\otimes\infty}+\tfrac12\mathrm{Ber}(0.8)^{\otimes\infty},
  \]
  where $\mathrm{Ber}(\theta)$ is the Bernoulli distribution with mean $\theta$.
  As shown in Section~\ref{sec:methods}, the measure in~\eqref{eq:numeraire} is $\Ps_t=\frac12\mathrm{Ber}(0.3)^{\otimes t}+\frac12\mathrm{Ber}(0.7)^{\otimes t}$ for $1\le t\le1000$, where $P^{\otimes t}$ is the distribution of $t$ independent variables with distribution $P$.
  Since the stopping time $\tau=1+\ind(X_1=1)$ gives $\E[\Es_\tau]=30/29>1$ under $\mathrm{Ber}(0.7)^{\otimes\infty}$, the numéraires do not form an e-process.
  The supremum $\eta_T$ in~\eqref{eq:snell} is attained at $\theta\in\{0.3,0.7\}$.
  Since $\tau\wedge T$ is a stopping time, $\eta_T$ is nondecreasing in $T$.
  In decimal arithmetic with $250$ digits, backward induction gives $9/8-\eta_T\approx4.7\times10^{-7}$, $3.3\times10^{-22}$, and $1.2\times10^{-100}$ for $T=50$, $200$, and $1000$.
  As a consequence, $\Gamma_T\le\log\eta_{1000}<\log(9/8)$ for $T\le1000$, although $\E_Q[\log\Es_T]$ grows linearly in $T$.
  Moreover, $\eta_T$ appears to increase to $9/8$.
  For smaller horizons, Table~\ref{tab:bernoulli} gives the values of $\Gamma_T$ computed by convex optimization and, for $T\in\{3,6\}$ and the uniform $\xi\in\Delta_T$, the value of~\eqref{eq:avg_dual} (Section~\ref{sec:methods}).
\end{example}

\begin{table}
  \tbl{Bernoulli example of \citet[Section~4.5]{koolen2022log}: optimal e-power at the horizon $T$, the value of~\eqref{eq:avg_dual} for the uniform $\xi\in\Delta_T$, and $\Gamma_T$, in nats}
  {%
    \begin{tabular}{rccc}
      $T$ & $\E_Q[\log\Es_T]$ & $\inf_{(E_t)}\sum_{t=0}^T\xi_tR_t(E_t)$ & $\Gamma_T$\\
      3 & 0.0464 & 0.0145 & 0.0260\\
      6 & 0.1332 & 0.0328 & 0.0456\\
      8 & 0.1908 & & 0.0515\\
      10 & 0.2471 & & 0.0547\\
      12 & 0.3019 & & 0.0563
  \end{tabular}}
  \label{tab:bernoulli}
\end{table}

\subsection{Stopping at the second observation}\label{sec:notep}

In the Gaussian location model of Section~\ref{sec:gauss}, the following proposition gives the conditional expectation of $\Es_2$ given $X_1$ and the best stopping time up to $t=2$ (Section~\ref{sec:numeraires}).

\begin{proposition}\label{prop:notep}
  For every $\theta\in\R$ and every $s\neq1$,
  \begin{gather}
    \E_\theta[\Es_2\mid X_1]=\{2/(1+s^2)\}^{1/2}\exp\{-(1-s^2)(X_1-\theta)^2/(2+2s^2)\},\label{eq:notep-cond}\\
    \sup_{\tau\le2}\E_\theta[\Es_\tau]=\E_\theta[\max\{1,\E_\theta[\Es_2\mid X_1]\}]>1,\label{eq:notep-sup}
  \end{gather}
  where the supremum runs over the stopping times of $(\cF_t)$ bounded by $2$.
\end{proposition}

Since the supremum in~\eqref{eq:notep-sup} exceeds one, the numéraires are not an e-process with respect to $(\cF_t)$ for $s\neq1$.
For $s=1/2$, this supremum equals $1.1132068$.
The rule that continues at $t=1$ exactly when $|X_1-\theta|<0.885$ attains it.
This threshold is slightly less than one standard deviation.
\citet[Section~3.2 and Supplementary Section~S3.1]{perezortiz2024estatistics} give a similar example for the $t$-test.
In contrast, \citet[Appendix~B.1]{grunwald2026supermartingales} show that, for $s<1$, the likelihood ratio $\prod_{i=2}^t\ell(Y_i)$ is a supermartingale with respect to $(\cG_t)$ also under normal distributions with variance larger than one.
Moreover, the null distributions of $X^2$ are absolutely continuous with respect to $Q_2$.
For $s\neq1$, Proposition~\ref{prop:minimax:i} and~\eqref{eq:notep-sup} therefore give $\Gamma_2>0$, and Proposition~\ref{prop:basic} gives $\Gamma_T\ge\Gamma_2>0$ for every $T\ge2$.
In particular, the minimax regret is positive for $s>1$, although Proposition~\ref{prop:sge1} bounds it by $\log s$.

\begin{proof}[of Proposition~\ref{prop:notep}]
  Under $P_\theta$, conditionally on $X_1$, $Y_2=2^{-1/2}(X_2-X_1)$ is normal with mean $2^{-1/2}(\theta-X_1)$ and variance $1/2$.
  Since $\Es_2=\ell(Y_2)$,~\eqref{eq:S-gauss} with $v=1/2$ gives~\eqref{eq:notep-cond}, because
  \[
    s^2(1+b/2)=(1+s^2)/2,\quad b/(2+b)=(1-s^2)/(1+s^2).
  \]
  Every stopping time $\tau\le2$ of $(\cF_t)$ satisfies $\{\tau\le1\}\in\cF_1$.
  Since $\Es_0=\Es_1=1$,
  \[
    \E_\theta[\Es_\tau]=\E_\theta[\ind(\tau\le1)+\ind(\tau=2)\E_\theta[\Es_2\mid X_1]]\le\E_\theta[\max\{1,\E_\theta[\Es_2\mid X_1]\}],
  \]
  with equality for $\tau=1+\ind\{\E_\theta[\Es_2\mid X_1]>1\}$.
  The conditional expectation $\E_\theta[\Es_2\mid X_1]$ exceeds one when
  \[
    (X_1-\theta)^2<\frac{1+s^2}{1-s^2}\log\frac{2}{1+s^2}
  \]
  if $s<1$, and when the reverse strict inequality holds if $s>1$.
  Since each of these conditions holds with positive probability under $P_\theta$, the right-hand side of~\eqref{eq:notep-sup} exceeds one.
\end{proof}

\subsection{An initial enlargement}\label{sec:enlarge}

It follows from Lemma~\ref{lem:orth} that $X^t$ and $(X_1,Y_2,\ldots,Y_t)$ are functions of each other in the Gaussian location model.
Therefore, for $t\ge1$,
\begin{equation}\label{eq:enlarge}
  \cF_t=\cG_t\vee\sigma(X_1).
\end{equation}
From $t=1$ on, the data filtration is thus the initial enlargement~\citep{jacod1985grossissement} of $(\cG_t)$ by the first observation.
We write $I(\cdot;\cdot)$ for the mutual information.
The following proposition shows that, under each null distribution, the first observation is a function of the whole sequence $(Y_i)_{i\ge2}$.
The variables $Y_2,\ldots,Y_t$ nevertheless carry only $\frac12\log t$ nats of information about it.
Moreover, for processes adapted to $(\cG_t)$, validity under one null distribution implies validity under every null distribution.

\begin{proposition}\label{prop:enlarge}
  \begin{parts}
  \item\partlabel{prop:enlarge:i} For every $\theta\in\R$, $X_1=\theta-\sum_{i\ge2}Y_i/\{i(i-1)\}^{1/2}$ almost surely under $P_\theta$.
  \item\partlabel{prop:enlarge:ii} Let $\theta\in\R$ and $t\ge1$.
    Under $P_\theta$, the conditional distribution of $X_1$ given $\cG_t$ is
    \[
      \N\biggl[\theta-\sum_{i=2}^t\frac{Y_i}{\{i(i-1)\}^{1/2}},\frac1t\biggr],
    \]
    and $I(X_1;Y^t)=\frac12\log t$.
  \item\partlabel{prop:enlarge:iii} Let $(E_t)$ be a process with values in $[0,\infty]$ that is adapted to $(\cG_t)$.
    Then $\sup_\tau\E_\theta[E_\tau]$ is the same for every $\theta\in\R$, where $\tau$ ranges over the stopping times of $(\cF_t)$.
    In particular, $(E_t)$ is an e-process for $\cP$ if and only if it is an e-process for $\{P_\theta\}$ for some $\theta\in\R$.
  \end{parts}
\end{proposition}

By Propositions~\ref{prop:numeraire} and~\ref{prop:enlarge:iii}, for every $\theta\in\R$, $\sup_\tau\E_\theta[\Es_{\tau\wedge T}]$ equals $1$ when $\tau$ ranges over the stopping times of $(\cG_t)$ and equals $\eta_T$ when $\tau$ ranges over the stopping times of $(\cF_t)$.
The variable $X_1$ in~\eqref{eq:enlarge} is thus the only source of the difference between the two suprema.

\begin{proof}[of Proposition~\ref{prop:enlarge}]
  (i) Under $P_\theta$, Lemma~\ref{lem:orth:ii} gives
  \[
    \bar X_t=X_1+\sum_{i=2}^t\frac{Y_i}{\{i(i-1)\}^{1/2}},
  \]
  and $\bar X_t\to\theta$ almost surely by the strong law of large numbers.

  (ii) By Lemma~\ref{lem:orth:iii}, $\bar X_t\sim\N(\theta,1/t)$ is independent of $Y^t$ under $P_\theta$.
  It follows that $X_1=\bar X_t-\sum_{i=2}^tY_i/\{i(i-1)\}^{1/2}$ has the stated conditional distribution.
  The mutual information is the difference $\frac12\log t$ between the differential entropies of $\N(0,1)$ and $\N(0,1/t)$.

  (iii) Let $\theta,\theta'\in\R$.
  Under $P_\theta$, the sequence $X'$ with terms $X'_i=X_i+\theta'-\theta$ has distribution $P_{\theta'}$.
  Each $E_t$ is a function of $Y^t$.
  By Lemma~\ref{lem:orth:i}, $Y^t$ takes the same value at $X$ and at $X'$, and so does $E_t$.
  For a stopping time $\tau$ of $(\cF_t)$, the map $x\mapsto\tau(x_1+\theta'-\theta,x_2+\theta'-\theta,\ldots)$ is also a stopping time of $(\cF_t)$.
  The value of $(E_t)$ at this stopping time has the same expectation under $P_\theta$ as $E_\tau$ under $P_{\theta'}$.
  Since the shift by $\theta'-\theta$ is a bijection of the set of stopping times of $(\cF_t)$, $\sup_\tau\E_\theta[E_\tau]=\sup_\tau\E_{\theta'}[E_\tau]$.
\end{proof}

\subsection{Insufficiency of a simple Gaussian alternative}\label{sec:simple}

The alternative $Q$ of Section~\ref{sec:gauss} mixes the alternatives $Q_\mu=\N(\mu,s^2)^{\otimes\infty}$ over the location.
The following proposition shows that this mixture is needed for Proposition~\ref{prop:mixture}: with $Q_\mu$ in place of $Q$, $\Gamma_T$ is bounded for every $s>0$.

\begin{proposition}\label{prop:simple}
  Let $s>0$, $\mu\in\R$, and $T\ge0$, and let $Q_\mu$ replace $Q$ in~\eqref{eq:numeraire} and~\eqref{eq:price}.
  \begin{parts}
  \item\partlabel{prop:simple:i} If $s\le1$, then $\Es_t=\prod_{i=1}^t\ell(X_i-\mu)$ and $\E_{Q_\mu}[\log\Es_t]=t\kappa$ for $t\ge0$.
    The process $(\Es_t)_{t\ge0}$ is a test supermartingale with respect to $(\cF_t)$.
    Hence, $\Gamma_T=0$.
  \item\partlabel{prop:simple:ii} If $s>1$, then $\Es_t=\prod_{i=2}^t\ell(Y_i)$ for $t\ge1$, as in Proposition~\ref{prop:numeraire}.
    Consequently, $\Gamma_T\le\log\eta_T\le\log s$, where $\eta_T$ is defined in~\eqref{eq:snell}.
  \end{parts}
\end{proposition}

Part~\ref{prop:simple:i*} is an instance of \citet[Theorem~1]{grunwald2024optimal}, since the variance $s^2$ of the alternatives is at most the variance of the null distributions.
For $s\le1$, the measure $\Ps_t$ in~\eqref{eq:numeraire} is the single null distribution $P_{\mu,t}$ at every horizon $t$.
In contrast, the mixing distribution $\N\{0,\rho^2-(1-s^2)/t\}$ of $\Ps_t$ in Proposition~\ref{prop:numeraire} changes with $t$.
This mixing variance is nonnegative at every horizon $t\ge1$ if and only if $\rho^2\ge1-s^2$, the condition of Section~\ref{sec:gauss}.
With $\rho=0$, the alternative of Section~\ref{sec:gauss} is $Q_\mu$ with $\mu=0$.
The condition $\rho^2\ge1-s^2$ then reduces to $s\ge1$.

\begin{proof}[of Proposition~\ref{prop:simple}]
  (i) Under every $P_\theta$, the factors of $\prod_{i=1}^t\ell(X_i-\mu)$ are independent, and~\eqref{eq:S-gauss} with $v=1$ gives
  \[
    \E_\theta[\ell(X_i-\mu)]=\exp\{-(1-s^2)(\theta-\mu)^2/2\}\le1.
  \]
  The product $\prod_{i=1}^t\ell(X_i-\mu)=\dd Q_{\mu,t}/\dd P_{\mu,t}$ is therefore a test supermartingale with respect to $(\cF_t)$.
  In particular, it belongs to $\cE_t$.
  Since $P_\mu\in\cP$, every $E\in\cE_t$ satisfies
  \[
    \E_{Q_\mu}\biggl[\frac{E}{\prod_{i=1}^t\ell(X_i-\mu)}\biggr]=\E_\mu[E]\le1.
  \]
  The product is therefore the numéraire $\Es_t$ of $Q_\mu$.
  Since the $X_i-\mu$ are independent $\N(0,s^2)$ variables under $Q_\mu$, $\E_{Q_\mu}[\log\Es_t]=t\kappa$.
  Since the numéraires form an e-process, Remark~\ref{rem:zero} gives $\Gamma_T=0$.

  (ii) For $t\ge1$, let $\Ps_t=\int P_{\theta,t}\N\{\mu,(s^2-1)/t\}(\dd\theta)$.
  Equation~\eqref{eq:S-factor} with $v=s^2$ and $m=\mu$ and its integral with $v=1$ over this mixing distribution give the densities of $Q_{\mu,t}$ and $\Ps_t$,
  \[
    t^{-1/2}\varphi_{s^2/t}(\bar x_t-\mu)\prod_{i=2}^t\varphi_{s^2}(y_i),\quad t^{-1/2}\varphi_{s^2/t}(\bar x_t-\mu)\prod_{i=2}^t\varphi_1(y_i).
  \]
  The proof of Proposition~\ref{prop:numeraire} then applies with $\N\{\mu,(s^2-1)/t\}$ in place of $\N\{0,\rho^2-(1-s^2)/t\}$ and gives $\Es_t=\prod_{i=2}^t\ell(Y_i)$.
  Since the null and $\Es_0,\ldots,\Es_T$ coincide with those of Section~\ref{sec:gauss}, the quantity $\eta_T$ in~\eqref{eq:snell} is also the same.
  Propositions~\ref{prop:basic:iii} and~\ref{prop:sge1} then give $\Gamma_T\le\log\eta_T\le\log s$.
\end{proof}

If $s<1$ and $0<\rho^2<1-s^2$, the mixing variance $\rho^2-(1-s^2)/t$ is negative at the first horizons, so that Proposition~\ref{prop:numeraire} does not apply there.
The following proposition shows that the rate of Proposition~\ref{prop:mixture:i} nevertheless holds for every $\rho>0$.

\begin{proposition}\label{prop:smallrho}
  Let $s\in(0,1)$ and $\rho>0$, and let $Q=\int Q_\mu\N(0,\rho^2)(\dd\mu)$ as in Section~\ref{sec:gauss}, without the condition $\rho^2\ge1-s^2$.
  Then $\Gamma_T/\log\log T\to1-s^2$ and $\log\eta_T/\log\log T\to1-s^2$ as $T\to\infty$.
\end{proposition}

With $\rho=0$, in contrast, Proposition~\ref{prop:simple:i} gives $\Gamma_T=0$ for every $T$.
Among the alternatives $Q_\mu$ and their mixtures over the location, $\Gamma_T$ thus tends to infinity exactly for the mixtures with $s<1$ and $\rho>0$, however small $\rho$ is.

\begin{proof}[of Proposition~\ref{prop:smallrho}]
  Let $\Pi_0=1$ and $\Pi_t=\prod_{i=2}^t\ell(Y_i)$ for $t\ge1$, and let $\nu_t=\sup_{E\in\cE_t}\E_Q[\log E]$.
  For every $\rho\ge0$, the proof of Proposition~\ref{prop:numeraire} shows that $(\Pi_t)$ is a test martingale with respect to $(\cG_t)$ and that $\E_Q[\log\Pi_t]=(t-1)\kappa$ for $t\ge1$.
  Hence, $\Pi_t\in\cE_t$ and $\nu_t\ge\E_Q[\log\Pi_t]$.
  For every $\theta\in\R$, \eqref{eq:S-gibbs}, applied to $Q_t$, $P_{\theta,t}$, and $E\in\cE_t$, gives $\nu_t\le\KL(Q_t\|P_{\theta,t})<\infty$.
  Let $t_\rho$ be the smallest integer with $t_\rho\ge(1-s^2)/\rho^2$.
  For $t\ge t_\rho$, the mixing variance $\rho^2-(1-s^2)/t$ is nonnegative, and the proof of Proposition~\ref{prop:numeraire} gives $\Es_t=\Pi_t$ and $\nu_t=\E_Q[\log\Pi_t]$.
  The differences $d_t=\nu_t-\E_Q[\log\Pi_t]$ therefore satisfy $0\le d_t\le C$ for every $t\ge0$, where $C=\max_{t<t_\rho}d_t<\infty$.

  For the upper bound, let $\eta'_T$ be the value of the stopping problem in Lemma~\ref{lem:identity:ii}.
  For every $\rho\ge0$, the proof of Lemma~\ref{lem:identity:ii} gives $\eta'_T=\sup_{\theta\in\R}\sup_\tau\E_\theta[\Pi_{\tau\wedge T}]$, where $\tau$ ranges over the stopping times of $(\cF_t)$.
  The process $(E_t)_{t\ge0}$ with $E_t=\Pi_t\ind(t\le T)/\eta'_T$ is therefore an e-process.
  By Lemma~\ref{lem:regret}, its regret at $t\le T$ is $d_t+\log\eta'_T\le C+\log\eta'_T$.
  Hence, $\Gamma_T\le C+\log\eta'_T$.
  The proof of the upper bound in Proposition~\ref{prop:mixture:i} gives $\limsup_{T\to\infty}\log\eta'_T/\log\log T\le1-s^2$, and hence $\limsup_{T\to\infty}\Gamma_T/\log\log T\le1-s^2$.

  For the lower bound, the proof of Lemma~\ref{lem:lbmix} holds with $\Pi_{t_k}$ in place of $\Es_{t_k}$, for every $\rho>0$.
  Since $\nu_{t_k}\ge\E_Q[\log\Pi_{t_k}]$, Lemma~\ref{lem:duality} still gives~\eqref{eq:S-lb}.
  Lemmas~\ref{lem:jk} and~\ref{lem:fixedk} and the proof of the lower bound in Proposition~\ref{prop:mixture:i} use only $\rho>0$ and~\eqref{eq:S-lb}, and they give $\liminf_{T\to\infty}\Gamma_T/\log\log T\ge1-s^2$.

  For the limit of $\log\eta_T$, Proposition~\ref{prop:basic:iii} gives $\Gamma_T\le\log\eta_T$, and the lower bound follows from that of $\Gamma_T$.
  For the upper bound, the distributions $P_{\theta,t}$ and $Q_t$ have positive densities on $\R^t$, and they therefore have the same null sets.
  For $t\ge t_\rho$, $\Es_t=\Pi_t$ $Q$-almost surely, and hence $P_\theta$-almost surely.
  For $t<t_\rho$, $\E_\theta[\Es_t]\le1$.
  Let $\tau$ be a stopping time of $(\cF_t)$.
  Since $\tau\vee t_\rho$ is also a stopping time of $(\cF_t)$,
  \[
    \E_\theta[\Es_{\tau\wedge T}]\le\sum_{t<t_\rho}\E_\theta[\Es_t]+\E_\theta[\Pi_{(\tau\vee t_\rho)\wedge T}]\le t_\rho+\eta'_T.
  \]
  Taking the supremum over $\theta$ and $\tau$ gives $\eta_T\le t_\rho+\eta'_T$.
  Since $\eta'_T\ge1$, it follows that $\eta_T\le(1+t_\rho)\eta'_T$, and hence $\limsup_{T\to\infty}\log\eta_T/\log\log T\le1-s^2$.
\end{proof}

\subsection{A test supermartingale with the running mean}

In each factor $\ell(X_i-\theta)$ of the likelihood ratio $L^\theta_t$ of Section~\ref{sec:stopping}, we replace the unknown mean $\theta$ by the running mean $\bar X_{i-1}$ and obtain the process $(\tilde E_t)$ of the following proposition.
Each factor $\ell(X_i-\bar X_{i-1})$ is the likelihood ratio of $\N(\bar X_{i-1},s^2)$ against $\N(\bar X_{i-1},1)$ at $X_i$.
\citet[Theorem~2, Part~3]{hao2026evalues} study a version of $(\tilde E_t)$ with a regularized running mean and obtain a regret of order $\{(1-s^2)/2\}\log t$ for each location.
For $(\tilde E_t)$ itself, the following proposition shows that the regret is an explicit harmonic sum.

\begin{proposition}\label{prop:plugin}
  Let $s<1$.
  The process $(\tilde E_t)_{t\ge0}$ with $\tilde E_t=\prod_{i=2}^t\ell(X_i-\bar X_{i-1})$ is a test supermartingale with respect to $(\cF_t)$ under every $P_\theta$, and, for $t\ge1$ and $T\ge1$,
  \[
    R_t(\tilde E_t)=\{(1-s^2)/2\}\sum_{i=1}^{t-1}i^{-1},\quad\Gamma_T\le\{(1-s^2)/2\}(1+\log T).
  \]
  For every $\mu$ and $t\ge1$, moreover, $R^\mu_t(\tilde E_t)=\kappa+R_t(\tilde E_t)$, with $R^\mu_t$ as in Theorem~\ref{thm:composite}.
\end{proposition}

\begin{proof}[of Proposition~\ref{prop:plugin}]
  Under $P_\theta$, conditionally on $\cF_{i-1}$, $X_i-\bar X_{i-1}$ is normal with mean $\theta-\bar X_{i-1}$ and variance $1$.
  Thus,~\eqref{eq:S-gauss} with $v=1$ gives
  \[
    \E_\theta[\ell(X_i-\bar X_{i-1})\mid\cF_{i-1}]=\exp\{-(1-s^2)(\theta-\bar X_{i-1})^2/2\}\le1.
  \]
  The supermartingale property follows.

  For the regret, $X_i-\bar X_{i-1}=\{i/(i-1)\}^{1/2}Y_i$, and
  \[
    \log\ell(y)-\log\ell[\{i/(i-1)\}^{1/2}y]=by^2/(2i-2).
  \]
  Since $\E_Q[Y_i^2]=s^2$ and $bs^2=1-s^2$, it follows that
  \[
    R_t(\tilde E_t)=\E_Q[\log\Es_t-\log\tilde E_t]=\sum_{i=2}^t\frac{1-s^2}{2(i-1)}.
  \]
  Since $\sum_{i=1}^{T-1}i^{-1}\le1+\log T$, the bound on $\Gamma_T$ follows.

  The distributions of $(Y_i)_{i\ge2}$ under $Q_\mu$ and under $Q$ coincide, so that the same computation holds under $Q_\mu$.
  Together with the e-power $t\kappa$ of the numéraire of $Q_\mu$ (Proposition~\ref{prop:simple:i}), it gives $R^\mu_t(\tilde E_t)=\kappa+R_t(\tilde E_t)$.
\end{proof}

\section{Numerical results}\label{sec:numerics}

Table~\ref{tab:snell} gives three upper bounds on $\Gamma_T$ for $T\le10^5$: an upper bound on $\log\eta_T$ (Proposition~\ref{prop:basic:iii} and Section~\ref{sec:methods}), $R_T(\tilde E_T)$ of Proposition~\ref{prop:plugin}, and $\log(T+1)$ of Proposition~\ref{prop:basic:ii}.
Up to $T=10^3$, the lower and upper bounds on $\log\eta_T$ of Section~\ref{sec:methods} coincide to five decimals.
At $T=10^4$ and $T=10^5$, they differ by at most $0.051$.
The exponentials of the upper bounds may replace $\eta_T$ in the e-process of Proposition~\ref{prop:basic:iii}.
At $T=10^5$, the upper bounds are at most $0.58$ times $(1-s^2)\log\log T$, the asymptotic rate of $\log\eta_T$ in Proposition~\ref{prop:mixture:i}.
At every horizon of Table~\ref{tab:snell}, $R_T(\tilde E_T)$ is more than three times the upper bound on $\log\eta_T$.

Table~\ref{tab:lb} gives Monte Carlo estimates of the limit as $\rho\to\infty$ of the right-hand side of~\eqref{eq:S-lb}, for the horizons $t_k=4^k$ with $k\le K$ and $4^K\le10^5$.
By Theorem~\ref{thm:composite:i}, these values estimate lower bounds on the minimax regret in the worst case over $\mu$ minus $\kappa$, and the upper bounds of Table~\ref{tab:snell} bound this quantity from above.
The estimates increase with $K$ and are positive at $K=8$ for the three values of $s$.
With the ratio $r=4$ of the horizons $t_k$, the lower bound $1-s^2(1-r^{-1/2})^{-2}$ on $\liminf_{T\to\infty}\Gamma_T/\log\log T$ in the proof of Proposition~\ref{prop:mixture:i} equals $1-4s^2$.

Since the numéraires form a test martingale with respect to $(\cG_t)$, we apply the adjuster $A(e)=(e-1-\log e)/(\log e)^2$ of \citet[Section~2.3]{choe2026combining} to their running maximum and obtain an e-process with respect to $(\cF_t)$.
We computed its regret for $s=1/2$ (Section~\ref{sec:methods}).
Since $\log A(e)=\log e-2\log\log e+o(1)$ as $e\to\infty$ and $\log\Es_t$ is close to $(t-1)\kappa$, we expect a regret close to $2\log\{(t-1)\kappa\}$.
At $t=10$, $10^2$, $10^3$, and $10^4$, the regret of the resulting e-process is $2.18$, $6.57$, $11.18$, and $15.80$, within 0.34 of $2\log\{(t-1)\kappa\}$.
The regret thus grows like $2\log t$, much faster than the regret of order $\log\log t$ of the e-process of Proposition~\ref{prop:mixture:ii}.

\begin{table}
  \tbl{Gaussian location model: upper bounds on $\Gamma_T$, in nats, rounded up: an upper bound on $\log\eta_T$, $R_T(\tilde E_T)$, and $\log(T+1)$}
  {%
    \begin{tabular}{rccccccc}
      & \multicolumn{2}{c}{$s=0.3$} & \multicolumn{2}{c}{$s=0.5$} & \multicolumn{2}{c}{$s=0.7$} & \\
      $T$ & $\log\eta_T$ & $R_T(\tilde E_T)$ & $\log\eta_T$ & $R_T(\tilde E_T)$ & $\log\eta_T$ & $R_T(\tilde E_T)$ & $\log(T+1)$\\
      $10$ & 0.3997 & 1.29 & 0.3195 & 1.07 & 0.2065 & 0.73 & 2.40\\
      $10^2$ & 0.6876 & 2.36 & 0.5550 & 1.95 & 0.3602 & 1.33 & 4.62\\
      $10^3$ & 0.9108 & 3.41 & 0.7390 & 2.81 & 0.4807 & 1.91 & 6.91\\
      $10^4$ & 1.1172 & 4.46 & 0.9004 & 3.68 & 0.5846 & 2.50 & 9.22\\
      $10^5$ & 1.2975 & 5.51 & 1.0425 & 4.54 & 0.6767 & 3.09 & 11.52
  \end{tabular}}
  \label{tab:snell}
\end{table}

\begin{table}
  \tbl{Gaussian location model: Monte Carlo estimates, in nats, of the limit as $\rho\to\infty$ of the lower bound in~\eqref{eq:S-lb}, for $T=4^K$ and $t_k=4^k$, with standard errors in parentheses}
  {%
    \begin{tabular}{rrccc}
      $K$ & $T$ & $s=0.3$ & $s=0.5$ & $s=0.7$\\
      1 & $4^{1}$ & $-0.128$ ($0.000$) & $-0.133$ ($0.001$) & $-0.133$ ($0.001$)\\
      2 & $4^{2}$ & $-0.017$ ($0.001$) & $-0.065$ ($0.001$) & $-0.099$ ($0.001$)\\
      3 & $4^{3}$ & $0.079$ ($0.001$) & $-0.012$ ($0.002$) & $-0.079$ ($0.002$)\\
      4 & $4^{4}$ & $0.180$ ($0.001$) & $0.056$ ($0.002$) & $-0.048$ ($0.003$)\\
      8 & $4^{8}$ & $0.514$ ($0.002$) & $0.247$ ($0.003$) & $0.030$ ($0.004$)
  \end{tabular}}
  \label{tab:lb}
\end{table}

\section{Numerical methods}\label{sec:methods}

For the Bernoulli example, let $\Ps_t$ be the measure of Example~\ref{ex:bernoulli}, and let $1\le t\le1000$.
For $\theta\in\{0.3,0.7\}$, the map $x\mapsto(1-x_1,\ldots,1-x_t)$ exchanges the two expectations $\E_\theta[\dd Q_t/\dd\Ps_t]$.
Since they average to $\E_{\Ps_t}[\dd Q_t/\dd\Ps_t]=1$, both equal one.
For $\theta\in\{0.4,0.5,0.6\}$, summation over the number of ones in decimal arithmetic with $60$ digits gives $\E_\theta[\dd Q_t/\dd\Ps_t]\le1$, with equality only at $t=1$.
Since $\Ps_t$ is also a mixture of null distributions, it is the measure in~\eqref{eq:numeraire}.
For each null parameter $\theta$, we define the Snell envelope of a process $(E_t)$ adapted to $(\cF_t)$ under $\mathrm{Ber}(\theta)^{\otimes\infty}$ by the recursion of Section~\ref{sec:setting}, with $(E_t)$ in place of $(\Es_t)$.
As in Section~\ref{sec:setting}, it is the smallest $\theta$-supermartingale that dominates $(E_t)$ on the binary tree of depth $T$~\citep[Theorem~1.2]{peskir2006optimal}.
Therefore, $(E_t)$ is an e-process for the five null distributions if and only if, for each $\theta$, some $\theta$-supermartingale $(\Psi^\theta_t)$ with $\Psi^\theta_t\ge E_t$ on the tree has initial value $\Psi^\theta_0\le1$.
With CVXPY and the Clarabel solver, we minimized $\max_{t\le T}R_t(E_t)$ over $(E_t)$ and the processes $(\Psi^\theta_t)$.
For the uniform $\xi$, we also minimized $\sum_t\xi_tR_t(E_t)$.
We computed the optimal e-powers $\E_Q[\log\Es_t]$ by exact summation.
For the uniform $\xi$, we reproduced the optimal values $0.002413$ and $0.021367$ of $\sup_{(E_t)}\sum_t\xi_t\E_Q[\log E_t]$ given by \citet[Section~4.5]{koolen2022log} to all six decimal places.
With Lemma~\ref{lem:regret}, these values give the third column of Table~\ref{tab:bernoulli}.
We computed the suprema $\eta_T$ by backward induction on the number of ones in decimal arithmetic with $250$ digits, separately for each of the five null parameters.

In the Gaussian location model, we computed lower and upper bounds on $\log\eta_T$ by backward induction for the stopping problem of Lemma~\ref{lem:identity:ii}, on a grid of times.
The grid contains every integer up to $1/0.0003$, and each later grid time is the integer part of $1.0003$ times the previous one.
Since the increment of $S_t$ between two grid times is normal, we integrated the transition of $S_t/t^{1/2}$ against the normal density on $3001$ equally spaced points of $[-12,12]$, with the value function divided by $\exp\{(1-s^2)S_t^2/(2t)\}$.
For the lower bound, the rule may stop only at grid times.
For the upper bound, the rule receives the largest reward of each block of grid times, and we bound this reward by the reflection principle for the Brownian bridge.
For $T\le10^3$, every block has length one, and the two bounds coincide.
For $s=0.3$ and $T=10^5$, a grid of $6001$ points on $[-16,16]$ leaves both bounds unchanged to five decimals.

For Table~\ref{tab:lb}, we computed the limit as $\rho\to\infty$ of the right-hand side of~\eqref{eq:S-lb} for the mixture $\bar P$ of Lemma~\ref{lem:lbmix} with $r=4$.
In this limit, $J_k$ of~\eqref{eq:S-J} tends to $\{p_k(Y^{t_k})+\xi_{t_k}\}/2$.
Under the alternative, we simulated the increments $\sum_{i=t_{k-1}+1}^{t_k}Y_i/\{i(i-1)\}^{1/2}$ directly as independent $\N\{0,s^2(1/t_{k-1}-1/t_k)\}$ variables.
We computed $p_k(Y^{t_k})$ by the midpoint rule in $z$ with $400$ points on each of the intervals $[-c-1,-c]$ and $[c,c+1]$.
After selecting $c$ on a grid of step $0.25$ with $1000$ simulated paths, we estimated the bound with $4000$ independent paths.
Fig.~\ref{fig:bounds} shows the same estimates.

For the comparison with adjusters in Section~\ref{sec:discussion}, we simulated $4000$ paths of $(\Es_t)_{t\le10^4}$ under $Q$ with $s=1/2$, with $Y_i\sim\N(0,s^2)$, and computed the average of $\log\Es_t-\log A(\max_{i\le t}\Es_i)$, with $A(e)=(e-1-\log e)/(\log e)^2$ for $e>1$ and $A(1)=1/2$.
The standard errors are at most $0.02$.

\section{Formalization in Lean}\label{sec:formalization}

\subsection{Overview}

We have formalized the results of the paper and of this Supplementary Material in the Lean proof assistant, with the Mathlib library.
We used Claude (Anthropic) to write the formalization, with our proofs as input.
The formalization is available at \url{https://adriencorenflos.github.io/formalization/log-optimal-e-processes}.
Every lemma, proposition, and theorem is stated and proved in Lean, as are the mathematical claims of the text that do not depend on a computation, except those of Section~\ref{sec:formalization-claims}.
We take the following results from the literature and assume them as hypotheses: the optimal stopping theorem of \citet[Theorem~1.2]{peskir2006optimal}, the primal-dual equivalence of \citet[Theorem~5]{koolen2022log}, the law of the iterated logarithm, the existence of numéraires~\citep[Theorem~2.6]{larsson2025numeraire}, and the e-lifting theorem of \citet[Theorem~2]{choe2026combining}.
We have formalized their statements but not their proofs.

We explain in Section~\ref{sec:formalization-improvements} how the formalization led us to improve three parts of the paper.
Section~\ref{sec:formalization-differences} describes the steps in which the formal proofs still differ from ours.

\subsection{Improvements of the paper}\label{sec:formalization-improvements}

Prior to the formalization, we proved Proposition~\ref{prop:minimax:ii} by applying Sion's minimax theorem to the sum $\sum_t\xi_tR_t(E_t)$ as a function of $\xi\in\Delta_T$ and of the e-process $(E_t)$ stopped at $T$.
This argument requires a topology on the e-processes in which the sum is lower semicontinuous.
When $\cX$ is finite, an e-process stopped at $T$ is a vector indexed by the nodes of the tree of depth $T$.
The sum is then continuous in this vector on the processes with terms $(1-\varepsilon)E_t+\varepsilon$.
Since we had no topology with this property on a general sample space, we stated Proposition~\ref{prop:minimax:ii} for finite $\cX$ only.
For the formalization, we instead applied Sion's theorem to the bilinear function $(\xi,v)\mapsto\sum_t\xi_tv_t$ on $\Delta_T$ and on the convex set of the vectors $v\in\R^{T+1}$ that dominate the regrets of some e-process.
Since this function is continuous, no topology on the e-processes is needed.
Proposition~\ref{prop:minimax:ii} now holds on every sample space.

Section~\ref{sec:rate} cites Proposition~\ref{prop:smallrho} for the claim that both limits of Proposition~\ref{prop:mixture:i} hold for every $\rho>0$.
Before the formalization, however, Proposition~\ref{prop:smallrho} gave only the limit of $\Gamma_T/\log\log T$.
When $\rho^2<1-s^2$, Proposition~\ref{prop:numeraire} does not give the numéraires of $Q$ at the horizons $t<(1-s^2)/\rho^2$.
Lemma~\ref{lem:identity:ii} therefore does not give $\eta_T$ for these mixtures.
We noticed this gap when we stated the claim in Lean.
Proposition~\ref{prop:smallrho} now also gives the limit of $\log\eta_T/\log\log T$.
In its proof, we bound $\eta_T$ by the value of the stopping problem of Lemma~\ref{lem:identity:ii} plus the number of horizons $t<(1-s^2)/\rho^2$.

We also shortened three of our proofs after the formalization.
In the proof of Lemma~\ref{lem:duality}, we previously bounded $\E_Q[\log^+E_t]$ by $\E_Q[\log\Es_t]+2$, through the proof of Lemma~\ref{lem:regret}, to show that the integrals in~\eqref{eq:S-gibbs} are well defined.
Since the functions $H$ and $f$ of that proof satisfy $\int(H/f)\dd Q^\xi\le1$, the positive part of $\log(H/f)$ is $Q^\xi$-integrable without this bound.
In our earlier proof of Remark~\ref{rem:deterministic}, we decomposed $\KL(Q^\xi\|\bar P)$ over the times $t\le T$ and applied~\eqref{eq:S-gibbs} at each time.
We now apply~\eqref{eq:S-gibbs} only once, on the pairs $(t,x^t)$, to the function equal to $\Es_t$ for $t\le T$ and to zero otherwise.
Finally, in our earlier proof of Proposition~\ref{prop:minimax:i}, we treated finite sample spaces by a separate argument.
By compactness, some e-process attained the minimax regret $\Gamma_T=0$.
The uniqueness of the log-optimal e-variable then showed that this e-process equals the numéraires $Q$-almost surely.
For absolutely continuous null distributions, we took e-processes with regrets at most $1/n$ and extracted a subsequence that converged almost surely to the numéraires.
With regrets at most $2^{-n}$, the whole sequence converges almost surely.
After the e-processes are set to zero at the nodes of $Q$-probability zero, the same argument applies to finite sample spaces.
We now use it in both cases.

\subsection{Steps in which the formal proofs differ}\label{sec:formalization-differences}

The other formal proofs follow ours, except in the following steps.
By Lemma~\ref{lem:orth:iii} and the convolution of normal distributions, the distributions of $(\bar X_t,Y^t)$ under the null distributions, their mixtures, and $Q$ are products of normal distributions.
For the formalization, we computed the densities of Proposition~\ref{prop:numeraire} and Lemma~\ref{lem:lbmix} from these products, instead of integrating~\eqref{eq:S-factor} over the mixing distributions.
In the proof of Proposition~\ref{prop:mixture:ii}, the strong law of large numbers gives $\Es_t\to0$ almost surely under $P_\theta$, while the argument requires only the almost sure convergence of $\Es_t$.
For the formalization, we instead obtained this convergence from Doob's martingale convergence theorem for the nonnegative martingale $(\Es_t)$ with respect to $(\cG_t)$.
The formalization also gives explicit terms $O(1)$ and $o(1)$ for both bounds in the proof of Proposition~\ref{prop:mixture:i}.
In the proof of Proposition~\ref{prop:enlarge:ii}, the mutual information is the difference between two differential entropies.
Since Mathlib has no differential entropy, we computed it in the formalization as the Kullback--Leibler divergence of the joint distribution of $(Y^t,X_1)$ from the product of its marginal distributions.

\subsection{Claims that are not formalized}\label{sec:formalization-claims}

Since the model of Example~\ref{ex:bernoulli} is not defined in the formalization, the symmetry argument for $\theta\in\{0.3,0.7\}$ in Section~\ref{sec:methods} is not formalized.
We have not formalized the identifications of our objects with those of the literature, such as the conditional e-variables of \citet{hao2026evalues} and the minimax pseudo-regret of \citet[Section~4]{grunwald2024safe}.

\end{document}